\documentclass[11pt,a4paper]{amsart}
\usepackage[utf8]{inputenc}
\usepackage[T1]{fontenc}
\usepackage{lmodern,microtype}
\usepackage{amsmath,amssymb,amsthm,mathtools}
\usepackage[margin=2.5cm]{geometry}
\usepackage{enumitem,booktabs,array,graphicx}
\usepackage{xcolor}
\usepackage[colorlinks=true,allcolors=black,bookmarksnumbered=true]{hyperref}
\hypersetup{pdftitle={Controllability of Acoustic Waves by Point Sources: Reachable States and Stable Reconstruction},pdfauthor={Kais Ammari and Mourad Sini}}
\setlist[itemize]{leftmargin=*,itemsep=3pt,topsep=5pt}
\setlist[enumerate]{leftmargin=*,itemsep=3pt,topsep=5pt}
\allowdisplaybreaks[2]
\newtheorem{theorem}{Theorem}[section]
\newtheorem{proposition}[theorem]{Proposition}
\newtheorem{lemma}[theorem]{Lemma}
\newtheorem{corollary}[theorem]{Corollary}
\theoremstyle{definition}

\newtheorem{remark}[theorem]{Remark}
\newtheorem{example}[theorem]{Example}

\numberwithin{equation}{section}

\newcommand{\R}{\mathbb R}
\newcommand{\C}{\mathbb C}
\newcommand{\Sph}{\mathbb S^2}
\newcommand{\dd}{\,\mathrm d}
\newcommand{\E}{\mathcal E}
\newcommand{\U}{\mathcal U_T}
\newcommand{\X}{\mathcal X}

\newcommand{\Reach}{\mathfrak R_T}
\newcommand{\Rad}{\mathsf R}
\newcommand{\Rr}{\mathcal R_T}
\newcommand{\RF}{\mathcal F_T}
\newcommand{\sinc}{\operatorname{sinc}}
\newcommand{\supp}{\operatorname{supp}}

\newcommand{\norm}[1]{\left\|#1\right\|}
\newcommand{\ip}[2]{\left\langle #1,#2\right\rangle}

\title[Point-source controllability and reconstruction]{Controllability of Acoustic Waves by Point Sources: Reachable States and Stable Reconstruction}
\author{Ka\"is Ammari}
\address{LR Analysis and Control of PDEs, LR22ES03, Department of Mathematics, Faculty of Sciences of Monastir, University of Monastir, Tunisia}
\email{kais.ammari@fsm.rnu.tn}
\author{Mourad Sini}
\address{RICAM, Austrian Academy of Sciences, Altenbergerstrasse 69, 4040 Linz, Austria}
\email{mourad.sini@oeaw.ac.at}
\date{11 November 2026}

\begin{document}
\begin{abstract}
We study the terminal controllability of the three-dimensional acoustic wave equation
\[
k^{-1}\partial_t^2w-\rho^{-1}\Delta w
=\sum_{j=1}^M q_j(t)\,\delta_{z_j}
\quad\text{in }(0,\infty)\times\R^3,
\qquad w(0)=\partial_tw(0)=0,
\]
where $k,\rho>0$, the points $z_j$ are distinct, and the amplitudes $q_j$ are supported in $[0,T]$. Although finitely many fixed sources cannot control the whole energy space, we establish the following precise results:
\begin{enumerate}[label=\textup{(\roman*)}]
\item For $q\in H^1(0,T)^M$, the terminal state has finite energy if and only if $q_j(0)=q_j(T)=0$ for every $j$.
\item The terminal energy satisfies the exact identity
\[
\E(S_Tq)=\frac{\kappa_0}{2\pi}\int_\R\omega^2
\big\langle G(\omega)\widehat q(\omega),\widehat q(\omega)\big\rangle\dd\omega,
\]
where $\kappa_0=\rho^2c/(4\pi k)$, $c=\sqrt{k/\rho}$, and
$G_{jj}=1$, $G_{j\ell}(\omega)=\sinc(\omega|z_j-z_\ell|/c)$ for $j\ne\ell$.
\item For every $T>0$, the map $S_T:H_0^1(0,T)^M\to\X=\dot H^1(\R^3)\times L^2(\R^3)$ is bounded below and is an isomorphism onto its closed reachable range $\Reach$. This range has infinite codimension even among states with the expected finite-propagation support.
\item We give an intrinsic radiation-field criterion for membership in
$\Reach$, reconstruct the unique controls stably, and solve the
best-approximation problem for non-reachable targets.
\end{enumerate}
The optimal coercivity constant equals one before the shortest source travel time; at long times it is of order $T^{-2}$ for two sources and $T^{-4}$ for three equally spaced collinear sources, with additional upper bounds for more general configurations.
Reproducible Galerkin computations illustrate these rates.
\end{abstract}
\maketitle

% === BEGIN INTRODUCTION (from supplied modular source) ===
\section{Introduction}\label{sec:introduction}

\subsection{The propagation problem}
Let $z_1,\ldots,z_M$ be distinct points in $\R^3$, let $k,\rho>0$ denote the bulk modulus and mass density of a homogeneous background, and set $c=\sqrt{k/\rho}$. We consider
\begin{equation}\label{eq:main}
\begin{cases}
k^{-1}\partial_t^2w-\rho^{-1}\Delta w
 =\displaystyle\sum_{j=1}^M q_j(t)\delta_{z_j},&t>0,\ x\in\R^3,\\
w(0)=\partial_tw(0)=0.
\end{cases}
\end{equation}
The source amplitudes $q_j$ are supported in a prescribed time interval $[0,T]$. Our objective is to characterize exactly the terminal states $(w(T),\partial_tw(T))$ that can be generated, to determine the topology in which their synthesis is stable, and to reconstruct the amplitudes from the terminal state.

Finite propagation already prevents exact or approximate controllability in the whole energy space at a fixed finite time. There is a further obstruction within the propagation region: finitely many fixed monopoles generate a restricted angular structure. These obstructions do not, however, describe the topology of the reachable set or the stability of source recovery. An injective propagation map could have a nonclosed range, in which case source synthesis would remain unstable even for reachable data. The issue addressed here is the stronger assertion that the finite-duration terminal map is bounded below in its natural amplitude norm.

The energy is
\begin{equation}\label{eq:energy-definition}
\E(u_0,u_1):=\frac12\int_{\R^3}
\left(\rho^{-1}|\nabla u_0|^2+k^{-1}|u_1|^2\right)\dd x.
\end{equation}
We equip
\[
\X=\dot H^1(\R^3)\times L^2(\R^3)
\quad\text{with}\quad \norm{y}_{\X}^2=\E(y),
\]
and use the amplitude space
\[
\U=H^1_0(0,T)^M,
\qquad \norm{q}_{\U}^2=\sum_{j=1}^M\norm{q_j'}_{L^2(0,T)}^2.
\]
The convention $\norm{y}_{\X}=\E(y)^{1/2}$ fixes all isometry constants below. The terminal map and its range are
\begin{equation}\label{eq:terminal-map}
S_Tq=(w(T),\partial_tw(T)),\qquad \Reach=S_T(\U).
\end{equation}
During emission the solution generally has a $1/|x-z_j|$ singularity and is not an energy-space trajectory. It is first constructed in a weaker space by the retarded potential. The endpoint conditions defining $\U$ are precisely those that produce finite terminal energy within the class of $H^1$ amplitudes. Exact reachability in this paper always refers to this terminal map; it does not assert surjectivity onto $\X$.

\subsection{Motivations}

\subsubsection{Resonant bubbles and the distinction between propagation and actuation}
One motivation comes from the time-dependent scattering of acoustic waves by small resonant bubbles. Mukherjee and Sini~\cite{MukherjeeSiniJEE} derive a leading scattered field represented by retarded monopoles whose amplitudes obey an effective dynamical system. With a suitable diagonal normalization $q=\Lambda_bY$, the propagation is described by~\eqref{eq:main}, while the amplitudes satisfy equations of the schematic form
\begin{equation}\label{eq:motivation-system}
d_iY_i''(t)+Y_i(t)+\sum_{j\ne i}b_{ij}Y_j''(t-\tau_{ij})
 =F_i[u^{\rm in}](t),\qquad
\tau_{ij}=|z_i-z_j|/c.
\end{equation}
Here $d_i=\omega_{M,i}^{-2}>0$ is the inverse squared Minnaert frequency of the isolated $i$th bubble, determined by its geometry and material parameters, while $b_{ij}=\mathcal C_j/(4\pi|z_i-z_j|)$ ($i\ne j$) measures the retarded interaction from bubble $j$ in the usual monopole normalization, with $\mathcal C_j$ its effective scattering coefficient. The forcing is induced by the prescribed incident field: in this normalization, $F_i[u^{\rm in}](t)=\partial_t^2u^{\rm in}(z_i,t)$, the second time derivative of the incident pressure at the bubble center. For instance, exterior point transducers produce $u^{\rm in}(z_i,t)$ as a superposition of distance-weighted, time-delayed copies of their prescribed signals; the effective source histories $q_j$ then emerge from the coupled bubble response.

This description leads naturally to two levels of control. At the \emph{ideal level}, the source amplitudes $q_j$ are prescribed independently, and the question is to characterize and reconstruct the terminal states generated through the propagation map $S_T$. At the \emph{physical level}, the amplitudes are constrained by the effective bubble dynamics and must be generated by admissible incident waves. The present paper is devoted to the ideal level.  The tracking study by Mukherjee and Sini ~\cite{MukherjeeSiniTracking} constructs ideal amplitudes for finite-dimensional projected pressure trajectories and studies their realization by resonant microstructures; its whole-space projection retains the boundary contribution of the observation region. The optimal-control study by Mukherjee and Sini ~\cite{MukherjeeSiniOptimal} considers finite-band wave--actuator systems and the cost of the input driving the internal amplitudes. Our current terminal-state problem uses the complete outgoing field, without a finite spectral projection. In particular, uniqueness of the full-state source history is consistent with nonuniqueness after a finite-dimensional observation. Passing from ideal to physical control is an important further step: it requires characterizing the amplitudes that can be produced by incident fields, estimating the associated control costs, and relating the resulting terminal states to waves scattered by finite-size bubbles. This second level is not addressed here and will be the subject of forthcoming work.

\subsubsection{Sound reproduction}
Point-source models also arise in loudspeaker arrays, active noise control, and sound-field reproduction. In these applications one chooses source histories to cancel an unwanted field or to synthesize a prescribed acoustic field in an observation region; see Nelson and Elliott~\cite{NelsonElliott} for active sound control and Onofrei and Platt~\cite{OnofreiPlatt} for acoustic source synthesis with controlled near fields. The terminal-state problem considered here concerns the complete source-generated field at a fixed time. Its range characterization describes which targets are compatible with emission by a finite array of fixed monopoles, and its inverse estimate quantifies the stability of the required amplitude histories.

\subsection{Related results and the contribution of this paper}
\paragraph{{\it{Concentrated control and spectral methods.}}}
Pointwise control of hyperbolic equations has a substantial history. Lions~\cite{Lions1984,Lions88,Lions92} develops point controls, transposition and their functional setting; Triggiani~\cite{Triggiani1993} studies interior and boundary regularity for wave equations with interior point control. In one dimension, Fabre and Puel~\cite{FabrePuel} relate pointwise controllability to limits of distributed control. The position and motion of the control are central in Khapalov~\cite{Khapalov95,Khapalov96} and Castro~\cite{Castro13}. Related coupled models include the boundary-controlled string with an interior point mass studied by Hansen and Zuazua~\cite{HansenZuazua}, and the flexible multi-structures treated by D\'ager and Zuazua~\cite{DagerZuazua}. The reduction of control and observation problems to moment equations and nonharmonic Fourier series is developed systematically by Avdonin and Ivanov~\cite{AvdoninIvanov} and Komornik and Loreti~\cite{KomornikLoreti}. Pointwise and spectral control of plate vibrations is studied by Haraux and Jaffard~\cite{HarauxJaffard}. These results involve different equations, boundary conditions, spectral hypotheses and mechanisms of actuation. They illustrate why statements about point control must specify the state space and the observation or control geometry. The present problem is the finite-duration, fixed-source propagation problem in three-dimensional free space.
\bigskip

\paragraph{{\it{Pointwise stabilization and actuator placement}}.}
Concentrated feedback is another closely related line of research. Ammari and Tucsnak study pointwise stabilization of Bernoulli--Euler beams~\cite{AmmariTucsnak00} and stabilization of abstract second-order equations by unbounded feedbacks~\cite{AmmariTucsnak01}. Singular internal stabilization of the wave equation is investigated by Jaffard, Tucsnak and Zuazua~\cite{JaffardTucsnakZuazua}; a broader treatment of elastic systems with collocated feedback is given by Ammari and Nicaise~\cite{AmmariNicaise}. These works address admissibility, observability and decay under feedback, whereas the source histories in~\eqref{eq:main} are prescribed controls over a finite interval. The location of a concentrated actuator can strongly affect its performance. Ammari, Henrot and Tucsnak~\cite{AmmariHenrotTucsnak} relate actuator placement to the asymptotic behavior of a pointwise stabilized string, and Privat, Tr\'elat and Zuazua~\cite{PrivatTrelatZuazua13} study optimal locations of distributed controllers for the one-dimensional wave equation. Here the geometry enters the terminal energy through the pairwise travel times and the angular Gram matrix. This supplies a different, explicit cost functional for the placement question discussed in Section~\ref{sec:conclusion}.
\bigskip

\paragraph{{\it{Identification of point sources.}}}
Inverse point-source problems provide a complementary perspective. Bruckner and Yamamoto~\cite{BrucknerYamamoto} study determination from pointwise observations, while Komornik and Yamamoto~\cite{KomornikYamamoto02,KomornikYamamoto05} establish estimates for identifying point sources and their inverse-problem applications. El Badia and Ha-Duong~\cite{ElBadiaHaDuong} use boundary measurements for point wave sources in three dimensions, and Al Jebawy, El Badia and Triki~\cite{AlJebawyElBadiaTriki} consider a moving point source. For vibrating beams, Nicaise and Za\"ir~\cite{NicaiseZair} treat identifiability, stability and reconstruction from boundary measurements. In the present paper the source locations are known, and the data for full reconstruction are the terminal state or its complete radiation profile. Determining unknown locations, and eventually potential sources, from partial observations leads to additional uniqueness and stability questions.
\bigskip

\paragraph{{\it{Radiation representation and the finite-time contribution.}}}
We use the free-wave radiation representation of Lax and Phillips and Friedlander~\cite{LaxPhillips,Friedlander80}, with its normalization and invertibility proved below by the Radon and Fourier transforms. Belishev and Vakulenko~\cite{BelishevVakulenko2008} study reachable and unreachable subspaces for three-dimensional acoustic scattering in the presence of a potential, with incoming waves prescribed at infinity as controls; they relate the failure of controllability to the existence of reversing waves. Their control mechanism differs from the finitely many fixed interior point sources considered here. Background on the Radon transform and abstract control duality can be found in~\cite{Helgason,TucsnakWeiss}. The angular sinc interaction kernel is also familiar in acoustic radiation: the power of a harmonic monopole array is expressed through a matrix with this kernel; see, for example,~\cite{SoaresDebut2025}.

\bigskip

\noindent {\it{Main contributions.}} Our principal contribution in this work is the explicit characterization of the range of the finite-duration terminal map. In particular, this characterization shows that neither exact nor approximate controllability holds in the prescribed ambient space; see Theorem~\ref{thm:main}. In addition, Theorem~\ref{thm:charact} gives necessary and sufficient conditions on the radiation field of an arbitrary energy state; Proposition~\ref{prop:frequency-membership} separates angular compatibility from the finite-support and zero-mean conditions on the recovered source derivatives. The classical radiation representation provides coordinates for the target state, whereas these additional conditions characterize precisely the states generated by the fixed source array. Compact temporal support also yields a coercive energy estimate despite the degeneracy of $G(0)$, making reconstruction on this closed, infinite-codimensional range stable. Finally, projection follows from this structure, while the geometry-dependent estimates quantify its conditioning. The main contribution is thus the range characterization with its finite-time stability and reconstruction.

\subsection{Main theorem}
Put
\[
a_{j\ell}=|z_j-z_\ell|/c\quad(j\ne\ell),\qquad
\kappa_0=\frac{\rho^2c}{4\pi k}.
\]
For $M\ge2$ write $d_{\min}=\min_{j\ne\ell}|z_j-z_\ell|$; for $M=1$ use the convention $d_{\min}=+\infty$. Define the real symmetric matrix
\begin{equation}\label{eq:G-definition}
G_{jj}(\omega)=1,\qquad
G_{j\ell}(\omega)=\sinc(\omega a_{j\ell}),\qquad
\sinc x=\begin{cases}\sin x/x,&x\ne0,\\1,&x=0.\end{cases}
\end{equation}
For compactly supported energy data, let $\Rad$ be the spatial Radon transform and set
\begin{equation}\label{eq:intro-radiation}
\RF[u_0,u_1](\sigma,\theta)
 =\partial_\sigma\left[
 (c^{-1}\Rad u_1-\partial_p\Rad u_0)
 \bigl(c(T-\sigma),\theta\bigr)\right].
\end{equation}
This map extends to all of $\X$, as established in Section~\ref{sec:radiation}.

\begin{theorem}[Terminal reachability and synthesis]\label{thm:main}
Let $T>0$ and let $z_1,\ldots,z_M$ be distinct.
\begin{enumerate}[label=\textup{(\roman*)}]
\item For $q\in H^1(0,T)^M$, the retarded solution has finite energy at $T$ if and only if $q_j(0)=q_j(T)=0$ for every $j$.
\item The map $S_T:\U\to\X$ satisfies
\begin{equation}\label{eq:main-bounds}
\kappa_0c_T\norm{q}_{\U}^2
\le\E(S_Tq)\le M\kappa_0\norm{q}_{\U}^2,
\end{equation}
where the optimal constant $c_T\in(0,1]$ obeys
\begin{equation}\label{eq:main-cutoff}
c_T\ge\frac12\inf_{|\omega|\ge\pi/(2T)}\lambda_{\min}(G(\omega)) >0,
\end{equation}
where $\lambda_{\min}(G(\omega))$ denotes the smallest eigenvalue of $G(\omega)$.

\medskip

If $M\ge2$ and $cT\le d_{\min}$, then $c_T=1$ and $\E(S_Tq)=\kappa_0\norm{q}_{\U}^2$. For $M=1$ the latter identity holds for every $T$.
\item The range $\Reach$ is closed, and $S_T:\U\to\Reach$ is a bounded isomorphism. It has infinite codimension even in the energy subspace supported in $K_T=\bigcup_j\overline{B(z_j,cT)}$.
\item A state $y\in\X$ belongs to $\Reach$ if and only if there is $f\in L^2(0,T)^M$, extended by zero, such that
\begin{equation}\label{eq:main-characterization}
\RF[y](\sigma,\theta)
 =\rho\sum_{j=1}^Mf_j(\sigma+\theta\cdot z_j/c),
\qquad \int_0^Tf_j(s)\dd s=0.
\end{equation}
The unique amplitude is $q_j(t)=\int_0^t f_j(s)\dd s$.
\end{enumerate}
\end{theorem}

Parts (i)--(iii) are proved in Sections~\ref{sec:setting}--\ref{sec:range}; part (iv) is proved in Section~\ref{sec:radiation}. Section~\ref{sec:reconstruction} gives the reconstruction equations and the orthogonal projection of arbitrary targets. Section~\ref{sec:quantitative} studies the time and geometry dependence of $c_T$, including two-sided long-time estimates for a pair and a collinear triple.
% === END INTRODUCTION ===

% === BEGIN SETTING (from supplied modular source) ===
\section{Distributional solutions and terminal energy}\label{sec:setting}

\subsection{Conventions and the retarded potential}
All amplitudes on $(0,T)$ are extended by zero to $\R$. The same symbol denotes the extension. We use
\[
\widehat f(\omega)=\int_\R f(t)e^{-i\omega t}\dd t,
\qquad
\widehat u(\xi)=\int_{\R^3}u(x)e^{-ix\cdot\xi}\dd x,
\]
so that the Plancherel factors are $(2\pi)^{-1}$ and $(2\pi)^{-3}$, respectively. The state and control spaces are real Hilbert spaces; complex Fourier transforms are used in the proofs. Write $r_j=|x-z_j|$ and $R_0=\max_j|z_j|$.

The space $\dot H^1(\R^3)$ is the completion of $C_c^\infty(\R^3)$ in the gradient norm. By the Sobolev inequality it is realized as
\[
\dot H^1(\R^3)=\{u\in L^6(\R^3):\nabla u\in L^2(\R^3)^3\}.
\]
Its dual is denoted by $\dot H^{-1}$, with norm
\[
\norm{v}_{\dot H^{-1}}
 =\sup_{\norm{\nabla\phi}_{L^2}=1}|\ip{v}{\phi}|.
\]
Every compactly supported element of $\dot H^1$ lies in $H^1$. The point mass $\delta_z$ does not belong to $H^{-1}(\R^3)$; the equation is therefore interpreted distributionally rather than as a standard energy-space forcing problem.

\begin{proposition}\label{prop:repr}
For $q\in L^2(0,T)^M$, the unique causal distributional solution of~\eqref{eq:main} is
\begin{equation}\label{eq:repr}
w(t,x)=\frac{\rho}{4\pi}\sum_{j=1}^M
\frac{q_j(t-r_j/c)}{r_j}.
\end{equation}
It satisfies
\begin{equation}\label{eq:weak-continuity}
(w,\partial_tw)\in C\bigl([0,\infty);L^2(\R^3)\times\dot H^{-1}(\R^3)\bigr),
\end{equation}
with zero initial state. There is a constant $C=C(M,k,\rho)$ such that
\begin{equation}\label{eq:weak-bound}
\sup_{t\ge0}\left(\norm{w(t)}_{L^2}+\norm{\partial_tw(t)}_{\dot H^{-1}}\right)
\le C\norm{q}_{L^2(0,T)^M}.
\end{equation}
For $t\ge T$, the solution is free and its spatial support is contained in
\begin{equation}\label{eq:shell-support}
\bigcup_{j=1}^M\{x:c(t-T)\le |x-z_j|\le ct\}.
\end{equation}
\end{proposition}

\begin{proof}
The operator $L=k^{-1}\partial_t^2-\rho^{-1}\Delta$ has causal fundamental solution
\[
\mathsf E(t,x)=\frac{\rho}{4\pi|x|}\delta(t-|x|/c).
\]
Convolution with $\sum_jq_j\otimes\delta_{z_j}$ gives~\eqref{eq:repr}. The existence and uniqueness of the causal fundamental solution, and hence uniqueness in causal distributions, are standard for the constant-coefficient wave equation; see~\cite{Hormander}. Alternatively one may first use smooth time functions, apply the usual distributional Duhamel formula, and pass to the limit in the estimates below.

For a single term $W_j=\rho q_j(t-r_j/c)/(4\pi r_j)$, polar coordinates give
\begin{equation}\label{eq:weak-L2}
\norm{W_j(t)}_{L^2}^2
=\frac{\rho^2c}{4\pi}\int_{-\infty}^{t}|q_j(s)|^2\dd s.
\end{equation}
Let $e_j(x)=(x-z_j)/r_j$. For smooth causal amplitudes, direct differentiation, including the wavefront in the distributional sense, gives
\begin{equation}\label{eq:weak-derivative}
\partial_tW_j=-c\operatorname{div}(e_jW_j)+cW_j/r_j.
\end{equation}
There is no additional mass at $z_j$: the flux of $e_jW_j$ through a sphere of radius $r$ tends to zero as $r\downarrow0$ for smooth amplitudes. Hardy's inequality in three dimensions yields
\[
|\ip{W_j/r_j}{\phi}|
\le\norm{W_j}_{L^2}\norm{\phi/r_j}_{L^2}
\le2\norm{W_j}_{L^2}\norm{\nabla\phi}_{L^2}.
\]
Together with~\eqref{eq:weak-derivative}, this gives
\[
\norm{\partial_tW_j}_{\dot H^{-1}}
\le3c\norm{W_j}_{L^2}.
\]
Density extends the formula and the bound to $L^2$ amplitudes. Summing and applying Cauchy--Schwarz to the finite sum proves~\eqref{eq:weak-bound}.

Continuity of time translation in $L^2(\R)$, followed by restriction to $r_j\ge0$, proves continuity of $W_j(t)$ in $L^2(\R^3)$ through~\eqref{eq:weak-L2}. Formula~\eqref{eq:weak-derivative} then proves continuity of its time derivative in $\dot H^{-1}$. At $t=0$ both limits are zero. Finally $q_j(t-r_j/c)$ can be nonzero only when $0\le t-r_j/c\le T$, proving~\eqref{eq:shell-support}. The right-hand side of the wave equation vanishes for $t>T$.
\end{proof}

\subsection{The radial energy calculation}
\begin{lemma}\label{lem:radial}
Let $f\in H^1(0,\infty)$ be compactly supported and satisfy $f(0)=0$. Then $W(x)=f(|x|)/|x|$ belongs to $H^1(\R^3)$, and
\begin{equation}\label{eq:radial-identities}
\norm{W}_{L^2}^2=4\pi\norm{f}_{L^2(0,\infty)}^2,
\qquad
\norm{\nabla W}_{L^2}^2=4\pi\norm{f'}_{L^2(0,\infty)}^2.
\end{equation}
\end{lemma}
\begin{proof}
For $f\in C_c^\infty(0,\infty)$, the first identity is polar integration. For the second,
\[
\frac1{4\pi}\norm{\nabla W}_{L^2}^2
=\int_0^\infty\left|f'-\frac fr\right|^2\dd r
=\int_0^\infty\left(|f'|^2-\frac{(|f|^2)'}r+\frac{|f|^2}{r^2}\right)\dd r.
\]
Integration by parts gives
\[
\int_0^\infty\frac{(|f|^2)'}r\dd r
=\int_0^\infty\frac{|f|^2}{r^2}\dd r,
\]
and proves the identity. Approximate a general zero-trace $f$ by $C_c^\infty(0,\infty)$ in $H^1$. The identities show convergence of the associated radial functions in $H^1(\R^3)$, and the first identity identifies their limit with $f(r)/r$. The boundary cancellation can also be checked directly from
\[
\frac{|f(r)|^2}{r}\le\int_0^r|f'(s)|^2\dd s\longrightarrow0.
\]
\end{proof}

\begin{theorem}[Exact endpoint condition within $H^1$]\label{thm:regularity}
Let $q\in H^1(0,T)^M$. The terminal state from Proposition~\ref{prop:repr} belongs to $\X$ if and only if
\begin{equation}\label{eq:endpoint-condition}
q_j(0)=q_j(T)=0\qquad (j=1,\ldots,M).
\end{equation}
Under~\eqref{eq:endpoint-condition},
\[
(w,\partial_tw)\in C([T,\infty);\X),\qquad
\E(w(t),\partial_tw(t))=\E(S_Tq)\quad(t\ge T),
\]
and
\begin{equation}\label{eq:rough-energy-upper}
\E(S_Tq)\le M\kappa_0\norm{q}_{\U}^2.
\end{equation}
\end{theorem}

\begin{proof}
\emph{Sufficiency.} A function in $H^1(0,T)$ has a zero extension in $H^1(\R)$ precisely when both endpoint traces vanish. For $t\ge T$, put
\[
f_j(r)=\frac{\rho}{4\pi}q_j(t-r/c),\qquad r\ge0.
\]
It is compactly supported, belongs to $H^1(0,\infty)$, and has $f_j(0)=0$. Lemma~\ref{lem:radial} gives $W_j(t)\in H^1(\R^3)$. Similarly,
\[
\partial_tW_j(t,x)=\frac{\rho}{4\pi r_j}q_j'(t-r_j/c)
\]
belongs to $L^2$. Directly from~\eqref{eq:radial-identities},
\begin{align*}
\E(W_j(t),\partial_tW_j(t))
&=2\pi\left(\frac{\rho}{4\pi}\right)^2
\left(\rho^{-1}c^{-1}+k^{-1}c\right)\norm{q_j'}_{L^2}^2\\
&=\kappa_0\norm{q_j'}_{L^2}^2.
\end{align*}
The inequality $\norm{\sum_j y_j}_{\X}^2\le M\sum_j\norm{y_j}_{\X}^2$ proves~\eqref{eq:rough-energy-upper}. The radial identities and continuity of translations of $q_j'$ imply continuity in $\X$ for $t\ge T$. The solution is free there, so energy conservation follows either by the Fourier representation of the free group or by approximation with smooth free data.

\emph{Necessity at the initial endpoint.} Suppose $w(T)\in\dot H^1$. For each $j$, the sphere $\Sigma_j=\{r_j=cT\}$ has an open surface patch avoiding all other spheres $\Sigma_\ell$ and all centers. Indeed distinct equal-radius spheres cannot share an open patch, and there are finitely many of them. In a neighborhood of this patch all terms with $\ell\ne j$ are $H^1$: their radial arguments avoid their jump interfaces and their denominators are bounded away from zero. The $j$th term has a jump of size $\rho q_j(0)/(4\pi cT)$ across $\Sigma_j$. An $H^1$ function has equal traces from the two sides of a smooth interior surface, so $q_j(0)=0$. This holds for every $j$.

\emph{Necessity at the terminal endpoint.} Once all initial traces vanish, the zero extension across $0$ is $H^1$ locally. Consequently every term with $\ell\ne j$ is $H^1$ in a sufficiently small neighborhood of $z_j$, including the case where another wavefront passes through that neighborhood. The remaining term can be written there as
\[
\frac{\rho q_j(T)}{4\pi r_j}
+\frac{\rho}{4\pi r_j}\bigl(q_j(T-r_j/c)-q_j(T)\bigr).
\]
The second summand is locally $H^1$ by Lemma~\ref{lem:radial}, with a cutoff. Since $\nabla(1/r_j)\notin L^2_{\rm loc}$, the first summand has a finite energy only when $q_j(T)=0$.
\end{proof}

%\begin{remark}[Energy during emission]\label{rem:active-energy}
%Applying the necessity argument on $(0,t)$ shows that, for $q\in H^1(0,T)^M$ and $0<t\le T$, finite energy at time $t$ is equivalent to $q_j(0)=q_j(t)=0$ for every $j$. In particular a nonzero amplitude history in $\U$ does not generate an energy-space trajectory throughout the control interval. The global weak solution of Proposition~\ref{prop:repr} and its finite-energy terminal state are the appropriate objects. This also explains why no energy identity is obtained by integrating an undefined pointwise product $q_j(t)\partial_tw(t,z_j)$ during emission.
%\end{remark}
%% === END SETTING ===

% === BEGIN ENERGY (from supplied modular source) ===
\section{An exact energy identity}\label{sec:energy}

For $a>0$, define the averaging operator on zero-extended functions by
\begin{equation}\label{eq:average}
A_av=\frac1{2a}\mathbf1_{[-a,a]}*v,
\qquad
(A_av)'(s)=\frac{v(s+a)-v(s-a)}{2a}.
\end{equation}
For $u,v\in H^1_0(0,T)$ put
\begin{equation}\label{eq:beta}
\beta_a(u,v)
=\int_\R u'(s)(A_av)'(s)\dd s
=\frac1{2a}\int_\R u'(s)\bigl(v(s+a)-v(s-a)\bigr)\dd s.
\end{equation}
The form is symmetric because $A_a$ is self-adjoint on $L^2(\R)$ and commutes with differentiation. Young's inequality gives
\begin{equation}\label{eq:beta-bound}
|\beta_a(u,v)|\le\norm{u'}_{L^2}\norm{v'}_{L^2}.
\end{equation}

\begin{lemma}\label{lem:archimedes}
For $e\in\R^3\setminus\{0\}$ and any function integrable on $[-|e|,|e|]$,
\begin{equation}\label{eq:archimedes}
\int_{\Sph}g(\theta\cdot e)\dd\theta
=\frac{2\pi}{|e|}\int_{-|e|}^{|e|}g(\tau)\dd\tau.
\end{equation}
\end{lemma}
\begin{proof}
Choose $e/|e|$ as the polar axis and put $\mu=\cos\vartheta$. The area element is $\dd\phi\dd\mu$ and $\theta\cdot e=|e|\mu$. Integrating first in $\phi\in[0,2\pi)$ proves the identity.
\end{proof}

\begin{theorem}[Temporal and angular energy identities]\label{thm:energy}
For every $q\in\U$,
\begin{align}
\E(S_Tq)
&=\frac{c\rho^2}{16\pi^2k}
\int_{\Sph}\int_\R
\left|\sum_{j=1}^Mq_j'(s+\theta\cdot z_j/c)\right|^2\dd s\dd\theta,
\label{eq:angular-energy}\\
&=\kappa_0\left(
\sum_{j=1}^M\norm{q_j'}_{L^2}^2
+\sum_{j\ne\ell}\beta_{a_{j\ell}}(q_j,q_\ell)\right).
\label{eq:energy}
\end{align}
The sums over $j\ne\ell$ are ordered sums. The right-hand sides depend on the zero-extended histories, not on the choice of a later terminal time after their support.
\end{theorem}

\begin{proof}
\emph{Step 1: the outgoing shell.} First take $q\in C_c^\infty(0,T)^M$. For $t>T+2R_0/c$, the support in~\eqref{eq:shell-support} lies in the origin-centered shell
\[
c(t-T)-R_0\le r=|x|\le ct+R_0.
\]
Its width is bounded independently of $t$. Uniformly on this shell, writing $x=r\theta$,
\[
r_j=r-\theta\cdot z_j+O(r^{-1}),\qquad
\frac{x-z_j}{r_j}=\theta+O(r^{-1}).
\]
Taylor expansion of the smooth, zero-extended amplitudes in~\eqref{eq:repr} gives
\begin{align*}
\partial_tw(t,r\theta)
 &=\frac{\rho}{4\pi r}\Phi(t-r/c,\theta)+O(r^{-2}),\\
\nabla w(t,r\theta)
 &=-\frac{\rho}{4\pi cr}\Phi(t-r/c,\theta)\theta+O(r^{-2}),
\end{align*}
where
\[
\Phi(s,\theta)=\sum_jq_j'(s+\theta\cdot z_j/c).
\]
All remainder constants are uniform in $(r,\theta)$ on the shell for these fixed smooth amplitudes. Since $\rho^{-1}c^{-2}=k^{-1}$, the energy density is
\[
\left(\frac{\rho}{4\pi}\right)^2k^{-1}r^{-2}|\Phi(t-r/c,\theta)|^2+O(r^{-3}).
\]
Integrating the error over a shell of fixed width and radius comparable to $t$ gives $O(t^{-1})$. The function $\Phi$ is supported in $[-R_0/c,T+R_0/c]$ in its first variable. The change $s=t-r/c$, followed by energy conservation after $T$, therefore yields~\eqref{eq:angular-energy}.

\emph{Step 2: interaction terms.} Expanding the square in~\eqref{eq:angular-energy}, each diagonal term integrates to $4\pi\norm{q_j'}_{L^2}^2$. For $j\ne\ell$ define
\[
C_{j\ell}(\tau)=\int_\R q_j'(s)q_\ell'(s+\tau)\dd s.
\]
With $e=(z_\ell-z_j)/c$, Lemma~\ref{lem:archimedes} gives
\begin{align*}
&\int_{\Sph}\int_\R
q_j'(s+\theta\cdot z_j/c)q_\ell'(s+\theta\cdot z_\ell/c)\dd s\dd\theta\\
&\hspace{12mm}=\frac{2\pi}{a_{j\ell}}
\int_{-a_{j\ell}}^{a_{j\ell}}C_{j\ell}(\tau)\dd\tau
=4\pi\beta_{a_{j\ell}}(q_j,q_\ell).
\end{align*}
The last equality follows by integrating $q_\ell'$ over the delay interval. Since $4\pi\,c\rho^2/(16\pi^2k)=\kappa_0$, this proves~\eqref{eq:energy}.

\emph{Step 3: density.} The map $S_T$ is bounded by Theorem~\ref{thm:regularity}. The angular expression is continuous on $\U$, because the square of a sum is bounded by $M$ times the sum of its squares and time translation preserves $L^2$. The temporal interaction expression is continuous by~\eqref{eq:beta-bound}. Approximation by $C_c^\infty(0,T)^M$ proves both identities for all $q\in\U$.
\end{proof}

\begin{corollary}[Frequency identity]\label{cor:fourier}
With $G$ defined in~\eqref{eq:G-definition},
\begin{equation}\label{eq:fourier-energy}
\E(S_Tq)=\frac{\kappa_0}{2\pi}
\int_\R\omega^2\ip{G(\omega)\widehat q(\omega)}{\widehat q(\omega)}\dd\omega.
\end{equation}
Here the Hermitian pairing is $\ip{v}{w}=\sum_jv_j\overline{w_j}$. For every $\xi\in\C^M$,
\begin{equation}\label{eq:G-Gram}
\ip{G(\omega)\xi}{\xi}
=\frac1{4\pi}\int_{\Sph}
\left|\sum_j\xi_j e^{i\omega\theta\cdot z_j/c}\right|^2\dd\theta.
\end{equation}
\end{corollary}
\begin{proof}
The Fourier transform of the averaging kernel is $\sinc(\omega a)$. Since the zero extension of $q\in\U$ belongs to $H^1(\R)^M$, its derivative has transform $i\omega\widehat q$ without endpoint masses. Plancherel applied to~\eqref{eq:energy} proves~\eqref{eq:fourier-energy}. Applying~\eqref{eq:archimedes} to the exponential of $i\omega\theta\cdot(z_j-z_\ell)/c$ gives~\eqref{eq:G-Gram}.
\end{proof}

\begin{lemma}\label{lem:G}
The matrix $G$ is continuous and even in $\omega$, and
\begin{enumerate}[label=\textup{(\roman*)}]
\item $0\le G(\omega)\le MI$ and $G(0)=\mathbf1\mathbf1^\top$;
\item $G(\omega)$ is positive definite for every $\omega\ne0$;
\item if $M\ge2$ and $d_{\min}=\min_{j\ne\ell}|z_j-z_\ell|$, then
\begin{equation}\label{eq:G-high-frequency}
\lambda_{\min}G(\omega)
\ge1-\frac{(M-1)c}{|\omega|d_{\min}},\qquad\omega\ne0.
\end{equation}
In particular $G(\omega)\to I$ as $|\omega|\to\infty$.
\end{enumerate}
\end{lemma}
\begin{proof}
Continuity, evenness and the value at zero follow from the sinc formula. The Gram identity proves nonnegativity. Cauchy--Schwarz in its finite sum gives the upper bound $MI$. The off-diagonal row sums are at most $(M-1)c/(|\omega|d_{\min})$, so the elementary Gershgorin estimate proves~\eqref{eq:G-high-frequency}.

For strict positivity, by evenness it suffices to consider $\omega>0$. If the Gram integral vanishes, its continuous angular integrand vanishes everywhere. Replacing $\theta$ by $-\theta$ gives
\[
\sum_j\xi_j e^{-i\varkappa\theta\cdot z_j}=0,
\qquad \theta\in\Sph,\qquad\varkappa=\omega/c>0.
\]
Consider the outgoing Helmholtz field
\[
v(x)=\sum_j\xi_j\frac{e^{i\varkappa|x-z_j|}}{4\pi|x-z_j|}.
\]
Its far-field pattern is zero. Rellich's lemma implies that $v$ is zero outside a large ball, and unique continuation extends this vanishing to the connected punctured domain $\R^3\setminus\{z_1,\ldots,z_M\}$; see~\cite{ColtonKress}. Near $z_j$, all other terms are smooth while $v(x)=\xi_j/(4\pi r_j)+O(1)$. Thus each $\xi_j$ is zero.
\end{proof}

%\begin{remark}[What the cross terms measure]
%The interaction terms in~\eqref{eq:energy} are interference terms in the radiated energy. They do not postulate a dynamical interaction between independently controlled sources. Dynamical coupling of physical actuators, when present, occurs in their separate amplitude equations. This distinction is useful when applying the identity to resonant bubbles.
%\end{remark}
% === END ENERGY ===

% === BEGIN RANGE (from supplied modular source) ===
\section{Coercivity and the topology of the reachable space}\label{sec:range}

\subsection{A finite-time spectral estimate}
\begin{theorem}[Quantitative coercivity]\label{thm:coercive}
For $b>0$, set
\[
m(b)=\inf_{|\omega|\ge b}\lambda_{\min}G(\omega).
\]
Then $m(b)>0$. The optimal constant
\begin{equation}\label{eq:cT-definition}
c_T:=\inf_{0\ne q\in\U}
\frac{\E(S_Tq)}{\kappa_0\norm{q}_{\U}^2}
\end{equation}
satisfies $0<c_T\le1$ and
\begin{equation}\label{eq:coercive}
\kappa_0c_T\norm{q}_{\U}^2\le\E(S_Tq)
\le M\kappa_0\norm{q}_{\U}^2.
\end{equation}
More precisely,
\begin{equation}\label{eq:cutoff-bound}
c_T\ge\sup_{0<b<\pi/T}
\left(1-\frac{bT}{\pi}\right)m(b)
\ge\frac12m\!\left(\frac{\pi}{2T}\right)>0.
\end{equation}
If all inter-source delays are at least $T$, the exact identity
\begin{equation}\label{eq:short-isometry}
\E(S_Tq)=\kappa_0\norm{q}_{\U}^2
\end{equation}
holds. For a single source it holds without a restriction on $T$.
\end{theorem}

\begin{proof}
By Lemma~\ref{lem:G}, $\lambda_{\min}G$ is continuous and strictly positive on compact subsets of $\R\setminus\{0\}$. It tends to one at infinity. Its infimum over $|\omega|\ge b$ is therefore positive.

Let $f=q'$, with zero extension. Its support is contained in $[0,T]$ and its Fourier transform is $i\omega\widehat q$. Cauchy--Schwarz, applied componentwise and then summed, gives
\[
|\widehat f(\omega)|^2
\le T\norm{f}_{L^2(0,T)^M}^2.
\]
There is no factor of $M$ in this estimate. Hence
\[
\frac1{2\pi}\int_{|\omega|<b}|\widehat f(\omega)|^2\dd\omega
\le\frac{bT}{\pi}\norm{f}_{L^2}^2.
\]
By~\eqref{eq:fourier-energy} and Plancherel,
\begin{align*}
\frac{\E(S_Tq)}{\kappa_0}
&\ge\frac{m(b)}{2\pi}
\int_{|\omega|\ge b}|\widehat f(\omega)|^2\dd\omega\\
&\ge m(b)\left(1-\frac{bT}{\pi}\right)\norm{q}_{\U}^2.
\end{align*}
This proves~\eqref{eq:cutoff-bound}. The upper bound in~\eqref{eq:coercive} follows from $G\le MI$. Taking only one nonzero amplitude gives a quotient exactly equal to one, so $c_T\le1$.

Finally, if $a\ge T$, then for almost every $s\in(0,T)$ one has $v(s+a)=v(s-a)=0$ for zero-extended $v\in H^1_0(0,T)$. Thus $\beta_a(u,v)=0$. When every pairwise delay is at least $T$, all cross terms in~\eqref{eq:energy} vanish. A single source has no cross terms.
\end{proof}

\begin{remark}
For $M\ge2$, $\lambda_{\min}G(\omega)$ tends to zero as $\omega\to0$. Thus a uniform pointwise lower bound in frequency is unavailable. Estimate~\eqref{eq:cutoff-bound} uses the fact that a signal supported in an interval of fixed length cannot place all its $L^2$ Fourier mass in an arbitrarily small interval about zero. It is this finite-time constraint that makes the terminal inverse bounded.
\end{remark}

%\begin{remark}[A compactness proof of coercivity]\label{rem:compactness-coercivity}
%There is also a qualitative proof that separates the high-frequency estimate from compactness on the control interval. Suppose that $\norm{q^n}_{\U}=1$ and $\E(S_Tq^n)\to0$. By Lemma~\ref{lem:G}, choose $\Omega>0$ such that $G(\omega)\ge\tfrac12 I$ for $|\omega|\ge\Omega$. The Fourier energy identity gives
%\[
%\int_{|\omega|\ge\Omega}\omega^2|\widehat q^n(\omega)|^2\dd\omega
%\le\frac{4\pi}{\kappa_0}\E(S_Tq^n)\longrightarrow0.
%\]
%After extraction, $q^n\rightharpoonup q$ in $\U$ and $q^n\to q$ in $L^2(0,T)^M$ by compactness of the embedding. The bounded linear map $S_T$ preserves weak convergence, so weak lower semicontinuity of the squared energy norm yields $\E(S_Tq)=0$. Positivity of $G(\omega)$ for $\omega\ne0$ in~\eqref{eq:fourier-energy} implies $\widehat q=0$ almost everywhere, and hence $q=0$.

%For the low frequencies,
%\[
%|\widehat q^n(\omega)-\widehat q(\omega)|
%\le\sqrt T\,\norm{q^n-q}_{L^2(0,T)^M}
%\]
%uniformly in $\omega$. Thus
%$\int_{|\omega|<\Omega}\omega^2|\widehat q^n(\omega)|^2\dd\omega\to0$.
%Combining the two frequency regions contradicts
%\[
%\int_\R\omega^2|\widehat q^n(\omega)|^2\dd\omega
%=2\pi\norm{q^n}_{\U}^2=2\pi.
%\]
%This proves the existence of a positive coercivity constant. Estimate~\eqref{eq:cutoff-bound} additionally makes its dependence on the time interval and the Gram matrix explicit.
%\end{remark}

\subsection{Closed range and an alternative operator interpretation}
\begin{theorem}[Closed reachable space]\label{thm:exact}
The terminal map $S_T$ is injective and has closed range in $\X$. Its inverse on $\Reach$ satisfies
\begin{equation}\label{eq:inverse-bound}
\norm{S_T^{-1}y}_{\U}
\le(\kappa_0c_T)^{-1/2}\norm{y}_{\X}.
\end{equation}
In the regime~\eqref{eq:short-isometry}, the map $\kappa_0^{-1/2}S_T$ is an isometry onto $\Reach$.
\end{theorem}
\begin{proof}
If $S_Tq=0$, coercivity gives $q=0$. If $S_Tq^n$ converges in $\X$, then
\[
\norm{q^n-q^m}_{\U}^2
\le(\kappa_0c_T)^{-1}\norm{S_Tq^n-S_Tq^m}_{\X}^2.
\]
Thus $q^n$ converges to some $q\in\U$. Continuity of $S_T$ identifies the state limit with $S_Tq$, proving closedness. The remaining assertions follow directly from the energy bounds and the stated norm convention.
\end{proof}

For later use introduce the symmetric form
\begin{equation}\label{eq:form-b}
b(q,v):=\sum_{\ell=1}^M\int_0^Tq_\ell'v_\ell'\dd s
+\sum_{j\ne\ell}\beta_{a_{j\ell}}(q_j,v_\ell).
\end{equation}
Polarization of the energy identity yields
\begin{equation}\label{eq:polarized-energy}
\ip{S_Tq}{S_Tv}_{\E}=\kappa_0 b(q,v),
\end{equation}
where $\ip{\cdot}{\cdot}_{\E}$ is the inner product whose squared norm is~\eqref{eq:energy-definition}.

\begin{proposition}\label{prop:compact-interaction}
There is a compact self-adjoint operator $\mathcal K_T$ on $\U$ such that
\[
b(q,v)=\ip{(I+\mathcal K_T)q}{v}_{\U}.
\]
The operator $I+\mathcal K_T$ is positive and invertible.
\end{proposition}
\begin{proof}
Riesz representation defines the bounded self-adjoint interaction operator. Symmetry of $\beta_a$ and its alternate representation give
\[
|\beta_a(u,v)|=|\beta_a(v,u)|
\le a^{-1}\norm{v'}_{L^2}\norm{u}_{L^2(0,T)}.
\]
Consequently $\norm{\mathcal K_Tq}_{\U}\le C\norm{q}_{L^2(0,T)^M}$. Compactness of $H^1_0(0,T)\hookrightarrow L^2(0,T)$ proves compactness of $\mathcal K_T$. Positivity and invertibility follow from Theorem~\ref{thm:coercive}. Conversely, the energy identity and Lemma~\ref{lem:G} show that the kernel is zero; the Fredholm alternative then gives another proof of invertibility and coercivity on the fixed interval.
\end{proof}

\subsection{A rigorous obstruction inside the propagation region}
The following source--wave identity will be used to prove that the reachable space has infinite codimension.

\begin{lemma}[Source--wave identity]\label{lem:duality-smooth}
Let $\varphi$ be a free wave with smooth compactly supported terminal data $\varphi(T)=\varphi_0$, $\partial_t\varphi(T)=\varphi_1$. For $q\in\U$,
\begin{equation}\label{eq:duality}
\frac1k\left(
\int_{\R^3}\partial_tw(T)\varphi_0\dd x
-\int_{\R^3}w(T)\varphi_1\dd x\right)
=\sum_{j=1}^M\int_0^Tq_j(t)\varphi(t,z_j)\dd t.
\end{equation}
\end{lemma}
\begin{proof}
For smooth time amplitudes, pair $Lw$ with the smooth free wave and integrate in time. The distributional spatial integration by parts is legitimate because $w$ is locally integrable and the test wave is smooth; a spatial cutoff can be chosen to contain their supports on $[0,T]$. The spatial terms cancel since $L\varphi=0$, leaving the displayed endpoint expression. This argument does not evaluate $w$ or $\partial_tw$ at a source. Proposition~\ref{prop:repr} provides the time traces needed for this integration, and the initial term is zero. Finally approximate $q$ in $\U$ by smooth compactly supported amplitudes; the terminal energy bound and the bounded smooth point traces pass the identity to the limit.
\end{proof}

Let
\[
\X_{K_T}=\{(u_0,u_1)\in\X:\supp u_0\cup\supp u_1\subset K_T\}.
\]
This is a closed subspace of $\X$.

\begin{proposition}[Failure of full controllability and infinite codimension]\label{prop:negative}
For every $T>0$, $\Reach\subset\X_{K_T}$, and
\[
\overline{\Reach}^{\,\X}\ne\X.
\]
Moreover $\Reach$ has infinite codimension in $\X_{K_T}$.
\end{proposition}
\begin{proof}
The support inclusion follows from Proposition~\ref{prop:repr}, and passes to energy limits. A nonzero smooth state supported outside $K_T$ cannot belong to this closed support subspace, proving the first assertion.

For the stronger statement, choose an integer $\ell\ge1$ such that $2\ell+1>M-1$. The real spherical harmonics of degree $\ell$ form a vector space of dimension $2\ell+1$. Hence there is a nonzero harmonic $Y_\ell$ satisfying
\[
Y_\ell\left(\frac{z_j-z_1}{|z_j-z_1|}\right)=0,
\qquad j=2,\ldots,M.
\]
These are at most $M-1$ homogeneous linear conditions. Fix $0<r_-<r_+<cT$. For any real $f\in C_c^\infty((r_-,r_+))$, set
\[
\varphi_0(x)=f(|x-z_1|)Y_\ell\left(\frac{x-z_1}{|x-z_1|}\right),
\qquad\varphi_1=0,
\]
with the smooth zero extension near the center. The free wave with these terminal data remains in the same angular harmonic: this follows by separation of the Laplacian in spherical coordinates, or by its commutation with the orthogonal projection onto that harmonic. It thus has the form $v(t,r)Y_\ell(\theta)$ for $r>0$. At the other centers its trace vanishes by the chosen angular zeros. At $z_1$ the trace vanishes by smoothness and $\ell\ge1$. Therefore $\varphi(t,z_j)=0$ for all $j$ and $t$.

The state $y_f=(0,\varphi_0)$ belongs to $\X_{K_T}$. By Lemma~\ref{lem:duality-smooth},
\[
\ip{S_Tq}{y_f}_{\E}
=\frac1{2k}\int_{\R^3}\partial_tw(T)\varphi_0\dd x=0
\qquad(q\in\U).
\]
The map $f\mapsto y_f$ is injective, since $Y_\ell$ is not zero. Its infinite-dimensional image is contained in $\X_{K_T}\cap\Reach^\perp$. This proves infinite codimension.
\end{proof}

\begin{example}\label{ex:one-source}
For a single source at the origin, the reachable states are exactly
\[
u_0(x)=\frac{f(r)}r,\qquad
u_1(x)=-\frac{cf'(r)}r,
\qquad f\in H^1_0(0,cT),
\]
extended by zero for $r>cT$. Indeed the representation gives $f(r)=\rho q(T-r/c)/(4\pi)$, and the converse is obtained from
\[
q(s)=\frac{4\pi}{\rho}f(c(T-s)).
\]
Thus radial symmetry alone is insufficient: the terminal pair must satisfy the outgoing relation $ru_1+c\partial_r(ru_0)=0$. This example illustrates the additional compatibility imposed by emission from rest.
\end{example}

\begin{corollary}\label{cor:nonzero-initial}
Let $W_0(t)$ denote the free wave group on $\X$. For initial state $y^{\rm in}\in\X$, the reachable terminal states are exactly
\[
W_0(T)y^{\rm in}+\Reach.
\]
A prescribed terminal state $y^{\rm out}$ can be reached precisely when $y^{\rm out}-W_0(T)y^{\rm in}\in\Reach$. The control is then unique and satisfies~\eqref{eq:inverse-bound} with this difference as target.
\end{corollary}
\begin{proof}
Linearity splits the solution into the free evolution of the initial data and the causal source-generated solution. Apply Theorem~\ref{thm:exact} to the latter.
\end{proof}
% === END RANGE ===

% === BEGIN RADIATION (from supplied modular source) ===
\section{Radiation fields and an intrinsic reachability criterion}\label{sec:radiation}

\subsection{Normalization and the free-wave representation}
For a smooth compactly supported spatial function let
\[
\Rad f(p,\theta)=\int_{x\cdot\theta=p}f(x)\dd S_x,
\qquad(p,\theta)\in\R\times\Sph.
\]
For smooth compactly supported data $y=(u_0,u_1)$ given at time $T$, define the radiation profile and its derivative by
\begin{align}
\Rr[y](\sigma,\theta)
&=\left(c^{-1}\Rad u_1-\partial_p\Rad u_0\right)
\bigl(c(T-\sigma),\theta\bigr),\label{eq:profile-definition}\\
\RF[y]&=\partial_\sigma\Rr[y].\label{eq:radiation-definition}
\end{align}
The derivative is the quantity controlled by energy. For general energy data the undifferentiated profile need not be an $L^2$ function; all statements for arbitrary targets will use $\RF$.

\begin{theorem}[Radiation transform]\label{thm:radiation}
Let $u$ be the free wave with data $y=(u_0,u_1)$ at time $T$.
\begin{enumerate}[label=\textup{(\roman*)}]
\item For smooth compactly supported data,
\begin{equation}\label{eq:radiation-limit}
\lim_{r\to\infty}r\,u(\sigma+r/c,r\theta)
=\frac1{4\pi}\Rr[y](\sigma,\theta).
\end{equation}
\item The map $\RF$ extends uniquely to a bounded bijection from $\X$ onto the real space $L^2(\R\times\Sph)$, and
\begin{equation}\label{eq:radiation-isometry}
\E(y)=C_{\rm rad}\norm{\RF[y]}_{L^2(\R\times\Sph)}^2,
\qquad C_{\rm rad}=\frac{c}{16\pi^2k}.
\end{equation}
\end{enumerate}
\end{theorem}

\begin{proof}
\emph{Step 1: the far-field limit.} Write $t=T+\tau$. Kirchhoff's formula is
\[
u(T+\tau,x)=\partial_\tau\bigl(\tau M_{c\tau}u_0(x)\bigr)
+\tau M_{c\tau}u_1(x),
\]
where $M_R$ is the normalized spherical mean. Set $x=r\theta$ and $\tau=\sigma-T+r/c$. On the compact support of the data the sphere $|y-r\theta|=c\tau$ converges to the plane $y\cdot\theta=c(T-\sigma)$. More explicitly, writing $y=p\theta+\eta$, $\eta\perp\theta$, its branch near the support is
\[
p=r-\sqrt{(r+c(\sigma-T))^2-|\eta|^2}
=c(T-\sigma)+\frac{|\eta|^2}{2r}+O(r^{-2}).
\]
On bounded $(\sigma,\eta)$ sets this convergence holds with the derivatives needed to differentiate once in $\sigma$; the surface Jacobian tends to one. Thus
\[
r\tau M_{c\tau}u_1(r\theta)
\longrightarrow\frac1{4\pi c}\Rad u_1(c(T-\sigma),\theta).
\]
For the $u_0$ term, differentiating the limiting plane location contributes $-c\partial_p$, while the derivative of the scalar prefactor tends to zero. This proves~\eqref{eq:radiation-limit}.

\emph{Step 2: the isometry.} The Fourier slice identity is
\[
\widehat{\Rad f}(\eta,\theta)=\widehat f(\eta\theta),
\]
with the transform on the left taken in $p$. The map $(\eta,\theta)\mapsto\eta\theta$ covers $\R^3\setminus\{0\}$ twice, and its radial Jacobian is $\eta^2$. Plancherel therefore gives
\begin{align*}
\int_{\Sph}\int_\R|\partial_p\Rad f|^2\dd p\dd\theta
&=8\pi^2\norm{f}_{L^2}^2,\\
\int_{\Sph}\int_\R|\partial_p^2\Rad f|^2\dd p\dd\theta
&=8\pi^2\norm{\nabla f}_{L^2}^2.
\end{align*}
The cross term between $\partial_p\Rad u_1$ and $\partial_p^2\Rad u_0$ integrates to zero: the first is odd and the second even under $(p,\theta)\mapsto(-p,-\theta)$. Since
\[
\RF[y]=-\partial_p\Rad u_1+c\partial_p^2\Rad u_0
\quad\text{at }p=c(T-\sigma),
\]
changing variables gives
\[
\norm{\RF[y]}_{L^2}^2
=\frac{8\pi^2}{c}\norm{u_1}_{L^2}^2
+8\pi^2c\norm{\nabla u_0}_{L^2}^2.
\]
This is~\eqref{eq:radiation-isometry}. Density of smooth compactly supported data gives an isometric extension, up to $C_{\rm rad}^{1/2}$, whose range is closed.

\emph{Step 3: surjectivity.} We give an explicit inverse, which also fixes all orientation and speed factors. Temporal Fourier transformation of~\eqref{eq:radiation-definition} gives
\begin{equation}\label{eq:radiation-fourier}
\widehat{\RF[y]}(\omega,\theta)
=\frac{e^{-i\omega T}}{c^2}
\left(i\omega\widehat u_1(-\omega\theta/c)
-\omega^2\widehat u_0(-\omega\theta/c)\right).
\end{equation}
For an arbitrary $F\in L^2(\R\times\Sph)$ and $\xi\ne0$, set $r=|\xi|$ and
\[
B_+(\xi)=e^{icrT}\widehat F(cr,-\xi/r),\qquad
B_-(\xi)=e^{-icrT}\widehat F(-cr,\xi/r).
\]
Define
\begin{equation}\label{eq:radiation-inverse}
\widehat u_0(\xi)=-\frac{B_+(\xi)+B_-(\xi)}{2r^2},
\qquad
\widehat u_1(\xi)=\frac{c\,[B_+(\xi)-B_-(\xi)]}{2ir}.
\end{equation}
The parallelogram identity gives, for $\xi\ne0$,
\[
 |\xi|^2|\widehat u_0(\xi)|^2+c^{-2}|\widehat u_1(\xi)|^2
 =\frac{|B_+(\xi)|^2+|B_-(\xi)|^2}{2|\xi|^2}.
\]
Thus polar coordinates, the changes of variables $\omega=\pm cr$,
and temporal Plancherel yield
\begin{align*}
 &\frac1{(2\pi)^3}\int_{\R^3}
 \left(|\xi|^2|\widehat u_0(\xi)|^2
       +c^{-2}|\widehat u_1(\xi)|^2\right)\dd\xi\\
 &\quad=\frac1{2(2\pi)^3}\int_{\Sph}\int_0^\infty
       \left(|B_+(r\theta)|^2+|B_-(r\theta)|^2\right)\dd r\dd\theta\\
 &\quad=\frac1{2c(2\pi)^3}\int_{\Sph}\int_\R
       |\widehat F(\omega,\theta)|^2\dd\omega\dd\theta
 =\frac1{8\pi^2c}\norm F_{L^2(\R\times\Sph)}^2.
\end{align*}
In particular $|\xi|\widehat u_0,\widehat u_1\in L^2(\R^3)$.
Since $k=\rho c^2$, spatial Plancherel gives
\[
 \E(u_0,u_1)=\frac1{16\pi^2\rho c}\norm F_{L^2}^2
 =\frac{c}{16\pi^2k}\norm F_{L^2}^2,
\]
with exactly the constant in~\eqref{eq:radiation-isometry}.
The weighted Fourier function $\widehat u_0$ represents an element of $\dot H^1$: it is locally integrable near zero by Cauchy--Schwarz and $|\xi|^{-1}\in L^2_{\rm loc}(\R^3)$, and the homogeneous Sobolev realization selects its $L^6$ representative. Substitution in~\eqref{eq:radiation-fourier} gives $\RF[y]=F$. 
If $F$ is real, then
$\widehat F(-\omega,\theta)=\overline{\widehat F(\omega,\theta)}$, and hence
\[
 B_+(-\xi)=\overline{B_-(\xi)},\qquad
 B_-(-\xi)=\overline{B_+(\xi)}.
\]
Formula~\eqref{eq:radiation-inverse} therefore gives
$\widehat u_j(-\xi)=\overline{\widehat u_j(\xi)}$ for $j=0,1$.
The recovered data are real, completing the proof of surjectivity.

\end{proof}

\begin{remark}[Relation to the translation representation]
Theorem~\ref{thm:radiation} is the free three-dimensional radiation representation, with the retarded time and acoustic constants specified. It is consistent with the Lax--Phillips and Friedlander theories~\cite{LaxPhillips,Friedlander80} and the Radon-transform identities used above are standard~\cite{Helgason}. The explicit inverse is included to avoid any additional support or parity assumption on an arbitrary derivative radiation field.
\end{remark}

\subsection{Reachable radiation profiles}
\begin{theorem}[Characterization of $\Reach$]\label{thm:charact}
For $y\in\X$, the following are equivalent:
\begin{enumerate}[label=\textup{(\roman*)}]
\item $y\in\Reach$;
\item there is $q\in\U$ such that
\begin{equation}\label{eq:derivative-characterization}
\RF[y](\sigma,\theta)
=\rho\sum_jq_j'(\sigma+\theta\cdot z_j/c)
\quad\text{for a.e. }(\sigma,\theta);
\end{equation}
\item there is $f\in L^2(\R)^M$, supported in $[0,T]$ and satisfying $\int_0^Tf_j=0$, such that $\RF[y]=\rho\sum_jf_j(\sigma+\theta\cdot z_j/c)$.
\end{enumerate}
The functions $q$ and $f=q'$ are unique, with $q_j(t)=\int_0^tf_j(s)\dd s$. In particular the criterion does not require a separate compact-support hypothesis on $y$; support in $K_T$ follows from it.
\end{theorem}
\begin{proof}
For $q\in C_c^\infty(0,T)^M$, its terminal state is smooth and compactly supported, because its radial functions vanish near both $r_j=0$ and $r_j=cT$. The free continuation for $t\ge T$ is~\eqref{eq:repr}. Therefore
\[
\lim_{r\to\infty}r\,w(\sigma+r/c,r\theta)
=\frac{\rho}{4\pi}\sum_jq_j(\sigma+\theta\cdot z_j/c).
\]
Comparison with~\eqref{eq:radiation-limit} and differentiation in $\sigma$ give~\eqref{eq:derivative-characterization}. Both sides are continuous maps from $\U$ to $L^2(\R\times\Sph)$, so density proves the formula for all $q\in\U$.

Conversely, if~\eqref{eq:derivative-characterization} holds, then $y$ and $S_Tq$ have the same derivative radiation field. The isometry implies $y=S_Tq$. The equivalence between (ii) and (iii) is the elementary isometry $q\mapsto q'$ from $H^1_0(0,T)$ onto the zero-mean subspace of $L^2(0,T)$. Uniqueness follows from Theorem~\ref{thm:exact}.
\end{proof}

\begin{corollary}\label{cor:profile-characterization}
For compactly supported $y\in\X$, the profile~\eqref{eq:profile-definition} is defined by approximation as a compactly supported $H^1$ function of $\sigma$ with values in $L^2(\Sph)$. Such a state is reachable if and only if, for some $q\in\U$,
\begin{equation}\label{eq:profile-characterization}
\Rr[y](\sigma,\theta)=\rho\sum_jq_j(\sigma+\theta\cdot z_j/c).
\end{equation}
\end{corollary}
\begin{proof}
Choose a ball containing the support of $y$ and approximate by smooth data supported in a fixed larger ball. Their profiles have a common bounded $\sigma$-support. By~\eqref{eq:radiation-isometry} their derivatives are Cauchy in $L^2$, and the one-dimensional Poincare inequality on that support gives convergence of the profiles in $H^1$. The limit agrees with the distributional Radon formula. The derivative characterization then determines the profiles up to a function of $\theta$ independent of $\sigma$; their compact support makes that difference zero.
\end{proof}

\begin{remark}
For $M=1$, $z_1=0$ and $q\in\U$, polar integration over a plane gives
\[
\Rad u_0(p)=\frac{\rho}{2}\int_{|p|}^{\infty}q(T-r/c)\dd r,
\qquad
\Rad u_1(p)=\frac{\rho c}{2}q(T-|p|/c).
\]
Also $\partial_p\Rad u_0(p)=-\rho\operatorname{sgn}(p)q(T-|p|/c)/2$ for a.e. $p$. Thus $c^{-1}\Rad u_1-\partial_p\Rad u_0=\rho q(T-p/c)$: for $p<0$ the left side cancels, and the right side is zero by the zero extension. This verifies~\eqref{eq:profile-characterization} on both signs of $p$.
\end{remark}

\subsection{A frequency-domain membership test}
For $\omega\ne0$ define
\[
V_\omega a(\theta)=\sum_ja_j e^{i\omega\theta\cdot z_j/c},
\qquad a\in\C^M.
\]
Then $(4\pi)^{-1}V_\omega^*V_\omega=G(\omega)$. The angular orthogonal projection onto $\operatorname{Ran}V_\omega$ is
\begin{equation}\label{eq:angular-projector}
P_\omega=\frac1{4\pi}V_\omega G(\omega)^{-1}V_\omega^*.
\end{equation}
This is an angular projection at each frequency; the separate temporal support condition below must still be imposed.

\begin{proposition}\label{prop:frequency-membership}
Let $F=\RF[y]$ for $y\in\X$, and set for a.e. $\omega\ne0$
\begin{equation}\label{eq:frequency-coefficients}
d(\omega)=\frac1{4\pi\rho}G(\omega)^{-1}V_\omega^*\widehat F(\omega,\cdot).
\end{equation}
Then $y\in\Reach$ if and only if
\begin{enumerate}[label=\textup{(\alph*)}]
\item $P_\omega\widehat F(\omega,\cdot)=\widehat F(\omega,\cdot)$ for a.e. $\omega$;
\item $d$ is the Fourier transform of an $f\in L^2(\R)^M$ supported in $[0,T]$ with $\int_0^Tf_j=0$ for every $j$.
\end{enumerate}
The source history is then $q_j(t)=\int_0^tf_j(s)\dd s$.
\end{proposition}
\begin{proof}
Fourier transformation of~\eqref{eq:derivative-characterization} yields $\widehat F=\rho V_\omega\widehat{q'}$. Applying $V_\omega^*$ proves necessity and identifies $d=\widehat{q'}$. Conversely (a),~\eqref{eq:angular-projector}, and~\eqref{eq:frequency-coefficients} give $\widehat F=\rho V_\omega d$. Condition (b) and Theorem~\ref{thm:charact} now imply reachability. The point $\omega=0$ has no effect on the $L^2$ identities. The real symmetry of $F$ and the definitions give a real $f$ whenever condition (b) holds.
\end{proof}
% === END RADIATION ===

% === BEGIN RECONSTRUCTION (from supplied modular source) ===
\section{Stable reconstruction and optimal approximation}\label{sec:reconstruction}

\subsection{Focusing the radiation data}
For arbitrary $F\in L^2(\R\times\Sph)$, define the derivative beamformed signals by
\begin{equation}\label{eq:derivative-beam}
(\mathcal BF)_\ell(s)
=\frac1{4\pi\rho}\int_{\Sph}
F(s-\theta\cdot z_\ell/c,\theta)\dd\theta,
\qquad \ell=1,\ldots,M.
\end{equation}
These are well-defined $L^2(\R)$ functions. Indeed Cauchy--Schwarz in $\theta$ and translation in $s$ give
\begin{equation}\label{eq:beam-bound}
\norm{\mathcal BF}_{L^2(\R)^M}
\le\frac{\sqrt M}{\rho\sqrt{4\pi}}\norm F_{L^2(\R\times\Sph)}.
\end{equation}
When $F=\RF[y]$ for a compactly supported state, the undifferentiated signals are
\begin{equation}\label{eq:beam-definition}
h_\ell(s)=\frac1{4\pi\rho}\int_{\Sph}
\Rr[y](s-\theta\cdot z_\ell/c,\theta)\dd\theta,
\qquad h_\ell'=(\mathcal BF)_\ell.
\end{equation}

\begin{proposition}[Beamforming equations]\label{prop:beamforming}
For $y=S_Tq$ with $q\in\U$, its signals satisfy on $\R$
\begin{align}
h_\ell&=q_\ell+\sum_{j\ne\ell}A_{a_{\ell j}}q_j,
\label{eq:beam}\\
(\mathcal B\RF[y])_\ell
&=q_\ell'+\sum_{j\ne\ell}A_{a_{\ell j}}q_j'.
\label{eq:derivative-beam-system}
\end{align}
In particular, for a.e. $\omega\ne0$,
\begin{equation}\label{eq:frequency-inverse}
\widehat q(\omega)=G(\omega)^{-1}\widehat h(\omega),
\qquad
\widehat{q'}(\omega)=G(\omega)^{-1}
\widehat{\mathcal B\RF[y]}(\omega).
\end{equation}
\end{proposition}
\begin{proof}
Insert~\eqref{eq:profile-characterization} in~\eqref{eq:beam-definition}. The source at $z_\ell$ contributes $q_\ell(s)$, while the contribution of $z_j$ is
\[
\frac1{4\pi}\int_{\Sph}
q_j(s+\theta\cdot(z_j-z_\ell)/c)\dd\theta
=\frac1{2a_{\ell j}}\int_{-a_{\ell j}}^{a_{\ell j}}q_j(s+\tau)\dd\tau.
\]
This is $A_{a_{\ell j}}q_j(s)$. Differentiating in the distributional, equivalently $L^2$, sense gives the second equation. Fourier transformation and Lemma~\ref{lem:G} prove~\eqref{eq:frequency-inverse}.
\end{proof}

\begin{remark}[Physical meaning of beamforming]\label{rem:beam-interpretation}
The direction variable $\theta$ may be viewed as indexing an ideal array of far-field receivers. In the normalized radiation profile $\rho^{-1}\Rr[y]$, a source at $z_\ell$ contributes $q_\ell(\sigma+\theta\cdot z_\ell/c)$, so the shift in~\eqref{eq:beam-definition} compensates its direction-dependent travel time. After that shift, this source contributes the same signal in every direction and therefore adds coherently under angular averaging. For another source at $z_j$, the remaining delay is $\theta\cdot(z_j-z_\ell)/c$. Angular averaging distributes that delay uniformly over $[-a_{\ell j},a_{\ell j}]$, by Lemma~\ref{lem:archimedes}. Its contribution is therefore the moving average $A_{a_{\ell j}}q_j$. In frequency, $G(\omega)$ is the matrix describing these interactions between focused signals. Its off-diagonal entries tend to zero at high frequency, whereas at zero frequency the columns coincide. For $M\ge2$, the finite temporal support constraint in Theorem~\ref{thm:coercive} is what permits stable recovery despite this low-frequency degeneracy.
\end{remark}

\noindent The following variational formulation enforces the temporal support and endpoint conditions and remains meaningful for arbitrary data, including data that do not satisfy the range conditions of Proposition~\ref{prop:frequency-membership}.

\subsection{A coercive variational inverse}
\begin{theorem}[Reconstruction from arbitrary radiation data]\label{thm:reconstruction}
Let $F\in L^2(\R\times\Sph)$ be real. There is a unique $q_F\in\U$ satisfying
\begin{equation}\label{eq:normal}
b(q_F,v)=\sum_{\ell=1}^M\int_0^T
(\mathcal BF)_\ell(s)v_\ell'(s)\dd s
\qquad\text{for every }v\in\U.
\end{equation}
Let $y_F=\RF^{-1}F$. Then
\begin{equation}\label{eq:projection}
S_Tq_F=P_{\Reach}y_F,
\end{equation}
where $P_{\Reach}$ is the energy-orthogonal projection. In particular, exact reachable data give the unique exact control. For any two data sets $F,H$,
\begin{align}
\norm{q_F-q_H}_{\U}
&\le \frac1{\rho\sqrt{4\pi c_T}}\norm{F-H}_{L^2(\R\times\Sph)},
\label{eq:radiation-noise-stability}\\
\norm{q_F-q_H}_{\U}
&\le c_T^{-1}\norm{\mathcal B(F-H)}_{L^2(0,T)^M}.
\label{eq:beam-noise-stability}
\end{align}
\end{theorem}

\begin{proof}
The right-hand side of~\eqref{eq:normal} is continuous on $\U$ by~\eqref{eq:beam-bound}. The form $b$ is bounded and coercive by~\eqref{eq:polarized-energy} and Theorem~\ref{thm:coercive}. The Lax--Milgram theorem therefore gives a unique $q_F$.

For $v\in\U$, polarize the radiation isometry and use Theorem~\ref{thm:charact}:
\begin{align*}
\ip{S_Tv}{y_F}_{\E}
&=C_{\rm rad}\rho\int_{\Sph}\int_\R
\sum_\ell v_\ell'(\sigma+\theta\cdot z_\ell/c)F(\sigma,\theta)\dd\sigma\dd\theta\\
&=4\pi C_{\rm rad}\rho^2\sum_\ell\int_0^T
v_\ell'(s)(\mathcal BF)_\ell(s)\dd s\\
&=\kappa_0\sum_\ell\int_0^T
v_\ell'(s)(\mathcal BF)_\ell(s)\dd s.
\end{align*}
The integrals are justified by Cauchy--Schwarz. In combination with~\eqref{eq:polarized-energy}, equation~\eqref{eq:normal} is precisely
\[
\ip{S_Tq_F-y_F}{S_Tv}_{\E}=0\qquad(v\in\U).
\]
Since $\Reach$ is closed, this proves~\eqref{eq:projection}.

The difference of the reconstructed states is the projection of $y_F-y_H$. Consequently
\[
\kappa_0c_T\norm{q_F-q_H}_{\U}^2
\le\norm{P_{\Reach}(y_F-y_H)}_{\X}^2
\le C_{\rm rad}\norm{F-H}_{L^2}^2.
\]
Since $C_{\rm rad}/\kappa_0=(4\pi\rho^2)^{-1}$, this proves~\eqref{eq:radiation-noise-stability}. Subtracting the two variational equations and testing with $v=q_F-q_H$ gives~\eqref{eq:beam-noise-stability} directly.
\end{proof}

\begin{corollary}[Best approximation of any terminal target]\label{thm:bestapprox}
For every $y\in\X$, the problem
\begin{equation}\label{eq:best-approx}
\min_{q\in\U}\E(S_Tq-y)
\end{equation}
has a unique solution $q_y=q_{\RF[y]}$. Its state is $P_{\Reach}y$, and
\begin{equation}\label{eq:state-noise-stability}
\norm{q_y-q_{\widetilde y}}_{\U}
\le(\kappa_0c_T)^{-1/2}\norm{y-\widetilde y}_{\X}.
\end{equation}
Moreover $y\in\Reach$ if and only if $\E(S_Tq_y-y)=0$.
\end{corollary}
\begin{proof}
Projection onto a closed subspace minimizes the squared Hilbert distance. Injectivity of $S_T$ makes the amplitude minimizer unique. Estimate~\eqref{eq:state-noise-stability} follows from the contraction property of the projection and~\eqref{eq:inverse-bound}.
\end{proof}

\begin{remark}\label{rem:affine-residual}
For a compactly supported target let $h$ be~\eqref{eq:beam-definition}, and define
\[
r_\ell=q_{y,\ell}+\sum_{j\ne\ell}A_{a_{\ell j}}q_{y,j}-h_\ell.
\]
The normal equation says $\int_0^Tr_\ell'v'=0$ for all $v\in H^1_0(0,T)$. The derivatives of such test functions are exactly the zero-mean functions in $L^2(0,T)$, so $r_\ell'$ is constant and $r_\ell$ is affine on $(0,T)$. For a reachable target the full-line beamforming identity makes this residual zero. For a general target it need not vanish; removing the affine terms would incorrectly discard the endpoint constraints. For a noncompact target one can take any local primitive of $(\mathcal B\RF[y])_\ell$ and obtain the same derivative formulation.
\end{remark}

\begin{proposition}[Quadratic regularization]\label{prop:tikhonov}
For $\varepsilon>0$, let $q_{y,\varepsilon}$ minimize
\[
\E(S_Tq-y)+\varepsilon\kappa_0\norm q_{\U}^2.
\]
Then
\begin{equation}\label{eq:tikhonov-rate}
\norm{q_{y,\varepsilon}-q_y}_{\U}
\le\frac{\varepsilon}{c_T+\varepsilon}\norm{q_y}_{\U}.
\end{equation}
\end{proposition}
\begin{proof}
Let $B_T=I+\mathcal K_T$ be the positive self-adjoint operator representing $b$. The unregularized and regularized equations have the same right-hand side, with left sides $B_Tq_y$ and $(B_T+\varepsilon I)q_{y,\varepsilon}$. Subtraction gives
\[
q_{y,\varepsilon}-q_y
=-\varepsilon(B_T+\varepsilon I)^{-1}q_y.
\]
Since $B_T\ge c_TI$, the inverse has norm at most $(c_T+\varepsilon)^{-1}$, proving the estimate.
\end{proof}

\subsection{Explicit formulas under additional separation}
\begin{proposition}[Neumann inversion and direct spherical recovery]\label{prop:explicit-inversion}
Let $y=S_Tq\in\Reach$ and let $h$ be its beamformed profile.
\begin{enumerate}[label=\textup{(\roman*)}]
\item Let $P$ restrict functions to $(0,T)$ and define, with zero extensions,
\[
(Kq)_\ell=\sum_{j\ne\ell}A_{a_{\ell j}}q_j.
\]
If $\eta:=\max_\ell\sum_{j\ne\ell}T/(2a_{\ell j})<1$, then
\begin{equation}\label{eq:neumann}
q=\sum_{n=0}^{\infty}(-PK)^nPh
\quad\text{in }L^2(0,T)^M.
\end{equation}
The error after terms $0,\ldots,N$ is at most $\eta^{N+1}\norm{Ph}_{L^2}/(1-\eta)$.
\item If $2cT\le d_{\min}$, then for a.e. $s\in(0,T)$,
\begin{equation}\label{eq:spherical-recovery}
q_j(s)=\frac{4\pi c(T-s)}{\rho}\,
\mathcal M_{z_j,c(T-s)}u_0,
\end{equation}
where $\mathcal M_{z,r}$ is the normalized spherical mean and $y=(u_0,u_1)$.
\end{enumerate}
\end{proposition}
\begin{proof}
The restricted averaging operator has kernel $(2a)^{-1}\mathbf1_{|s-t|\le a}$ on $(0,T)^2$. Its $L^2$ norm is at most $\min\{1,T/(2a)\}$ by Schur's test. The symmetric matrix of these block bounds has operator norm at most its largest row sum, so $\norm{PK}\le\eta$. Restricting~\eqref{eq:beam} gives $(I+PK)q=Ph$, and the usual geometric-series estimate proves (i).

Under the separation in (ii), the other source contributions vanish in the open ball $B(z_j,cT)$. Thus $u_0=\rho q_j(T-r_j/c)/(4\pi r_j)$ there. Its spherical average has the same value for a.e. radius, which proves~\eqref{eq:spherical-recovery}.
\end{proof}

%\begin{remark}[Compatibility and noise]
%Formula~\eqref{eq:frequency-inverse} is an exact identity on compatible data. Since $G(\omega)^{-1}$ is unbounded near zero for $M\ge2$, it is not an unrestricted bounded Fourier multiplier on arbitrary noisy signals. The stable inverse in Theorem~\ref{thm:reconstruction} enforces $q\in H^1_0(0,T)^M$ through its variational formulation. Likewise the $L^2$ Neumann inverse in~\eqref{eq:neumann}, applied to arbitrary data, does not by itself enforce endpoint vanishing or $H^1$ regularity. The variational formulation is the appropriate reconstruction when terminal energy and general perturbations of the data are part of the specification.
%\end{remark}
% === END RECONSTRUCTION ===

% === BEGIN QUANTITATIVE (from supplied modular source) ===
\section{Geometry and time dependence of the control cost}\label{sec:quantitative}

We use throughout the \emph{optimal} coercivity constant
\begin{equation}\label{eq:optimal-constant}
c_T:=\inf_{0\ne q\in\U}
\frac{\E(S_Tq)}{\kappa_0\sum_{j=1}^M\norm{q_j'}_{L^2(0,T)}^2}.
\end{equation}
Theorem~\ref{thm:coercive} gives $c_T>0$ for every fixed configuration of distinct sources and every finite $T$. The inverse estimate for source synthesis involves $c_T^{-1/2}$. Its dependence on the time horizon and the geometry is therefore relevant even though exact synthesis on the reachable subspace holds for all $T>0$.

\subsection{Elementary properties and scaling}

\begin{proposition}\label{prop:cost-properties}
For a fixed source configuration, $0<c_T\le1$ and $T\mapsto c_T$ is nonincreasing. If $M=1$, then $c_T=1$ for every $T$. If $M\ge2$ and
\[
a_*:=\min_{j\ne\ell}\frac{|z_j-z_\ell|}{c},
\]
then $c_T=1$ for $T\le a_*$. The constant depends only on the dimensionless pairwise delays $|z_j-z_\ell|/(cT)$. In particular, it is invariant under rigid motions of the source configuration, and
\[
c_{\lambda T}(\lambda z_1,\ldots,\lambda z_M;c)
=c_T(z_1,\ldots,z_M;c),\qquad \lambda>0.
\]
\end{proposition}

\begin{proof}
Positivity was proved in Theorem~\ref{thm:coercive}. A control with only one nonzero component has energy $\kappa_0\norm{q_j'}_2^2$, so $c_T\le1$. The same identity proves the assertion for $M=1$.

Let $0<T_1<T_2$. Zero extension embeds $H_0^1(0,T_1)^M$ isometrically into $H_0^1(0,T_2)^M$. For an embedded control, the source has stopped by $T_1$, and its wave energy is conserved between $T_1$ and $T_2$. Thus its quotient in \eqref{eq:optimal-constant} is the same at both horizons. Taking the infimum over the larger control space gives $c_{T_2}\le c_{T_1}$.

For completeness, if $a\ge T$ and $s\in(0,T)$, then the zero extension of $v\in H_0^1(0,T)$ satisfies $v(s+a)=v(s-a)=0$ for almost every $s$. Hence $\beta_a(u,v)=0$. Applying the energy identity to every off-diagonal term proves $c_T=1$ when $T\le a_*$. This argument uses cancellation in the energy form; it does not require the balls $B(z_j,cT)$ to be disjoint.

Finally, make the change of variables $s=T\sigma$ in the energy quotient. Both the diagonal derivative norms and every $\beta_{a_{j\ell}}$ acquire the factor $T^{-1}$, while the delays become $a_{j\ell}/T$. The common factor cancels. The remaining assertions follow because the energy identity involves only pairwise distances.
\end{proof}

\subsection{Dipole cancellation and the two-source rate}

\begin{proposition}[A dipole upper bound]\label{prop:dipole-bound}
If $M\ge2$, then
\begin{equation}\label{eq:dipole-bound}
c_T\le\min\left\{1,\frac{2\pi^2}{3}\left(\frac{a_*}{T}\right)^2\right\}.
\end{equation}
In particular, $c_T\to0$ as $T\to\infty$.
\end{proposition}

\begin{proof}
Choose a pair with delay $a_*$ and use the amplitudes $g$ and $-g$ at these two points, setting all other amplitudes to zero. We take $g\in H_0^2(0,T)$, where $H_0^2$ denotes the closure of $C_c^\infty(0,T)$ in $H^2$; hence $g$ and $g'$ vanish at both endpoints, and the zero extension belongs to $H^2(\R)$. The energy identity gives
\begin{align*}
\E(S_Tq)
&=2\kappa_0\bigl(\norm{g'}_2^2-\beta_{a_*}(g,g)\bigr)\\
&=\frac{\kappa_0}{\pi}\int_\R
\omega^2\bigl(1-\sinc(\omega a_*)\bigr)|\widehat g(\omega)|^2\dd\omega.
\end{align*}
The identity
\[
1-\sinc x=\int_0^1(1-\cos(tx))\dd t\le\frac{x^2}{6}
\]
holds for every real $x$. Plancherel's theorem therefore yields
\[
\E(S_Tq)\le\frac{\kappa_0a_*^2}{3}\norm{g''}_2^2.
\]
Since $\sum_j\norm{q_j'}_2^2=2\norm{g'}_2^2$, we obtain
\begin{equation}\label{eq:dipole-test}
c_T\le\frac{a_*^2}{6}\frac{\norm{g''}_2^2}{\norm{g'}_2^2}.
\end{equation}
For
\[
g(t)=\frac{T}{2\pi}\left(1-\cos\frac{2\pi t}{T}\right),
\]
we have $g\in H_0^2(0,T)$, $g'(t)=\sin(2\pi t/T)$, and
\[
\frac{\norm{g''}_2^2}{\norm{g'}_2^2}=\frac{4\pi^2}{T^2}.
\]
Substituting this value into \eqref{eq:dipole-test} proves the second upper bound. The first follows from Proposition~\ref{prop:cost-properties}.
\end{proof}

\begin{theorem}[Two-source coercivity at long times]\label{thm:two-source-rate}
For $M=2$, write $a=|z_1-z_2|/c$. Then
\begin{equation}\label{eq:two-source-explicit}
\frac{19}{240}\min\left\{\frac{\pi^2a^2}{4T^2},1\right\}
\le c_T
\le\min\left\{1,\frac{2\pi^2a^2}{3T^2}\right\}.
\end{equation}
Consequently,
\[
c_T\asymp (a/T)^2\qquad\text{as }T/a\to\infty.
\]
The comparison constants in this statement are independent of $a$ and $T$.
\end{theorem}

\begin{proof}
The two eigenvalues of $G(\omega)$ are $1\pm\sinc(\omega a)$, so
\[
\lambda_{\min}G(\omega)=1-|\sinc(\omega a)|.
\]
We first prove the elementary bound
\begin{equation}\label{eq:sinc-lower}
1-|\sinc x|\ge\frac{19}{120}\min\{x^2,1\},\qquad x\in\R.
\end{equation}
By evenness, it suffices to consider $x\ge0$. For $0\le x\le1$, the alternating Taylor estimate for $\sin x$ gives
\[
1-\sinc x\ge\frac{x^2}{6}-\frac{x^4}{120}
\ge\frac{19x^2}{120}.
\]
On $[1,\pi]$, $\sinc x$ is nonnegative and nonincreasing. Indeed,
$\sin x-x\cos x$ has derivative $x\sin x\ge0$ on $[0,\pi]$ and vanishes at zero. It follows that
\[
1-|\sinc x|\ge1-\sin1\ge\frac{19}{120},\qquad 1\le x\le\pi.
\]
For $x\ge\pi$, use $|\sinc x|\le1/\pi$ and $1-1/\pi>19/120$. This proves \eqref{eq:sinc-lower}.

The quantitative spectral estimate in Theorem~\ref{thm:coercive} now gives
\begin{align*}
c_T
&\ge\frac12\inf_{|\omega|\ge\pi/(2T)}\lambda_{\min}G(\omega)\\
&\ge\frac{19}{240}\min\left\{\left(\frac{\pi a}{2T}\right)^2,1\right\}.
\end{align*}
The upper bound is Proposition~\ref{prop:dipole-bound}. These estimates imply the asserted two-sided order.
\end{proof}

\begin{corollary}[Regular simplices]\label{cor:simplex-rate}
Suppose $M\in\{2,3,4\}$ and all pairwise source distances equal $d>0$. Then, with $a=d/c$,
\[
c_T\asymp (a/T)^2\qquad\text{as }T/a\to\infty.
\]
\end{corollary}

\begin{proof}
Writing $s=\sinc(\omega a)$ and $J$ for the all-ones matrix, we have $G=(1-s)I+sJ$. Its eigenvalues are
\[
1-s\quad(M-1\text{ times}),\qquad 1+(M-1)s.
\]
Thus $\lambda_{\min}G(\omega)\sim (\omega a)^2/6$ at zero. Positivity away from zero follows from Lemma~\ref{lem:G}, while $G(\omega)\to I$ at infinity. By continuity and compactness on the intermediate frequency interval, there is $\gamma_M>0$ such that
\[
\lambda_{\min}G(\omega)\ge\gamma_M\min\{(\omega a)^2,1\}.
\]
The spectral lower bound of Theorem~\ref{thm:coercive} and the dipole upper bound prove the result.
\end{proof}

%\begin{remark}[Order and leading coefficient]
%Theorems above establish the long-time \emph{order} of $c_T$ for the stated geometries. They do not identify a sharp leading coefficient. In particular, the dipole test gives
%\[
%\limsup_{T/a\to\infty}(T/a)^2c_T\le\frac{2\pi^2}{3}
%\]
%5for two sources. Equality and the existence of the corresponding limit are not asserted here.
%\end{remark}

\subsection{A different rate for a collinear triple}

\begin{theorem}[Quadrupole cancellation]\label{thm:collinear-rate}
Let $M=3$ and, up to a rigid motion, let
\[
z_1=-de_1,\qquad z_2=0,\qquad z_3=de_1,
\qquad d>0.
\]
Put $a=d/c$. Then
\begin{equation}\label{eq:triple-upper}
c_T\le\min\left\{1,\frac{2\pi^2}{3}\left(\frac aT\right)^2,
\frac{32\pi^4}{15}\left(\frac aT\right)^4\right\},
\end{equation}
and there is a constant $C_3>0$, independent of $a$ and $T$, such that
\begin{equation}\label{eq:triple-lower}
c_T\ge C_3\min\{(a/T)^4,1\}.
\end{equation}
In particular,
\[
c_T\asymp(a/T)^4\qquad\text{as }T/a\to\infty.
\]
\end{theorem}

\begin{proof}
\emph{Upper bound.}
Choose $q=(g,-2g,g)$, with $g\in H_0^3(0,T)$. Its angular source factor is
\[
e^{-i\omega a\theta_1}-2+e^{i\omega a\theta_1}
=2\bigl(\cos(\omega a\theta_1)-1\bigr).
\]
Since $|2(\cos x-1)|\le x^2$ and
$(4\pi)^{-1}\int_{\Sph}\theta_1^4\dd\theta=1/5$, the Gram representation and Plancherel give
\[
\frac{\E(S_Tq)}{\kappa_0}
\le\frac{a^4}{5}\norm{g'''}_2^2.
\]
The denominator in \eqref{eq:optimal-constant} is $6\norm{g'}_2^2$. Therefore
\begin{equation}\label{eq:triple-test}
c_T\le\frac{a^4}{30}\frac{\norm{g'''}_2^2}{\norm{g'}_2^2}.
\end{equation}
Take $g(t)=g_0(t/T)$ with $g_0(s)=(1-\cos(2\pi s))^2$ on $(0,1)$. The first three boundary traces $g_0,g_0',g_0''$ vanish, so $g_0\in H_0^3(0,1)$. Orthogonality of the sine functions in
\[
g_0'=4\pi\sin(2\pi s)-2\pi\sin(4\pi s),\qquad
g_0'''=-16\pi^3\sin(2\pi s)+32\pi^3\sin(4\pi s)
\]
yields
\[
\norm{g_0'}_2^2=10\pi^2,\qquad
\norm{g_0'''}_2^2=640\pi^6.
\]
Consequently $\norm{g'''}_2^2/\norm{g'}_2^2=64\pi^4/T^4$, and \eqref{eq:triple-test} gives the last upper bound in \eqref{eq:triple-upper}. The other two bounds have already been proved.

\emph{Lower bound.}
Set $x=\omega a$, $s=\sinc x$, and $t=\sinc(2x)$. The Gram matrix is
\[
G_3(x)=\begin{pmatrix}1&s&t\\s&1&s\\t&s&1\end{pmatrix}.
\]
Reflection of the two end sources decomposes $\R^3$ into an antisymmetric line and a symmetric plane. On the line, the eigenvalue is $1-t$. In the orthonormal basis $2^{-1/2}(1,0,1)$, $(0,1,0)$ of the plane, the matrix is
\[
\begin{pmatrix}1+t&\sqrt2s\\\sqrt2s&1\end{pmatrix},
\]
with eigenvalues
\[
\lambda_\pm(x)=1+\frac t2\pm\frac12\sqrt{t^2+8s^2}.
\]
Taylor expansion gives
\[
s=1-\frac{x^2}{6}+\frac{x^4}{120}+O(x^6),\qquad
t=1-\frac{2x^2}{3}+\frac{2x^4}{15}+O(x^6).
\]
In particular, $1-t\sim2x^2/3$, $\lambda_+(x)\to3$, and
\[
\lambda_-(x)\lambda_+(x)
=1+t-2s^2
=\frac{2}{45}x^4+O(x^6).
\]
It follows that
\begin{equation}\label{eq:triple-small-frequency}
\lambda_{\min}G_3(x)=\lambda_-(x)
=\frac{2}{135}x^4+O(x^6)
\qquad(x\to0).
\end{equation}
By Lemma~\ref{lem:G}, $G_3(x)$ is positive definite for $x\ne0$, and it tends to $I$ as $|x|\to\infty$. Hence
\[
\gamma_3:=\inf_{x\ne0}
\frac{\lambda_{\min}G_3(x)}{\min\{x^4,1\}}>0.
\]
Indeed, \eqref{eq:triple-small-frequency} bounds this quotient below near zero, convergence to $I$ bounds it below for large $|x|$, and continuity gives a positive minimum on the remaining compact set. Theorem~\ref{thm:coercive} therefore implies
\[
c_T\ge\frac{\gamma_3}{2}
\min\left\{\left(\frac{\pi a}{2T}\right)^4,1\right\}
\ge\frac{\gamma_3}{2}\min\{(a/T)^4,1\}.
\]
This proves \eqref{eq:triple-lower} with $C_3=\gamma_3/2$.
\end{proof}

\noindent Theorems~\ref{thm:two-source-rate} and~\ref{thm:collinear-rate} show that the long-time exponent cannot be inferred from the minimum separation alone. Spatial moment cancellation changes the order at which a source combination becomes visible in the radiated field.

\begin{proposition}[A sufficient higher-order cancellation criterion]\label{prop:multipole-upper}
Let $p\ge1$, let $z_*\in\R^3$, and suppose that $v=(v_1,\ldots,v_M)\in\R^M\setminus\{0\}$ satisfies
\begin{equation}\label{eq:moment-cancellation}
\sum_{j=1}^M v_j(z_j-z_*)^\alpha=0
\quad\text{for every multi-index }\alpha\text{ with }|\alpha|<p.
\end{equation}
Then, for every nonzero $g\in C_c^\infty(0,1)$ and every $T>0$,
\begin{equation}\label{eq:multipole-upper}
c_T\le\frac{A_p^2}{|v|^2T^{2p}}
\frac{\norm{g^{(p+1)}}_{L^2(0,1)}^2}{\norm{g'}_{L^2(0,1)}^2},
\qquad
A_p:=\frac1{p!c^p}\sum_{j=1}^M|v_j|\,|z_j-z_*|^p.
\end{equation}
\end{proposition}

\begin{proof}
Use the control $q_j(t)=v_jg(t/T)$ and write $g_T(t)=g(t/T)$. The multinomial formula shows that \eqref{eq:moment-cancellation} cancels every term of order less than $p$ in
\[
\sum_jv_j\exp\bigl(i\omega\theta\cdot(z_j-z_*)/c\bigr).
\]
The integral remainder in Taylor's formula for $e^{ix}$ has modulus at most $|x|^p/p!$ for real $x$. Therefore the displayed sum has modulus at most $A_p|\omega|^p$, uniformly in $\theta\in\Sph$. Translation by $z_*$ contributes only a common phase and does not affect its modulus. The Gram identity now yields
\[
\frac{\E(S_Tq)}{\kappa_0}
\le A_p^2\norm{g_T^{(p+1)}}_2^2,
\qquad
\sum_j\norm{q_j'}_2^2=|v|^2\norm{g_T'}_2^2.
\]
The quotient of derivative norms scales as $T^{-2p}$, which proves \eqref{eq:multipole-upper}.
\end{proof}

\noindent Proposition~\ref{prop:multipole-upper} gives an upper bound by testing a fixed moment-cancelling source combination. A matching lower bound must instead control all combinations, including frequency-dependent ones. In particular, a spectral estimate
\[
 \lambda_{\min}G(\omega)\ge C\min\{|\omega|^{2p},1\}
\]
for a fixed configuration would give $c_T\gtrsim T^{-2p}$ for large $T$ by~\eqref{eq:cutoff-bound}. Establishing its optimal exponent for arbitrary configurations requires further small-frequency spectral analysis that is beyond the scope of the present paper. The pair, regular-simplex and collinear-triple results supply explicit two-sided examples; the general classification is separate from the reachable-space characterization.

\begin{corollary}[Attainment of the optimal coercivity constant]
\label{cor:cost-attainment}
For every $T>0$, the infimum in~\eqref{eq:optimal-constant} is attained.
If $c_T<1$, it is an eigenvalue of finite multiplicity of
$B_T=I+\mathcal K_T$ on $\U$: there exists $q_*\in\U$ such that
\[
 \norm{q_*}_{\U}=1,\qquad B_Tq_*=c_Tq_*,\qquad
 \E(S_Tq_*)=\kappa_0c_T.
\]
\end{corollary}
\begin{proof}
By Proposition~\ref{prop:compact-interaction}, $\mathcal K_T$ is compact and self-adjoint, and
\[
 c_T-1=\inf_{\norm q_{\U}=1}\ip{\mathcal K_Tq}{q}_{\U}
       =\inf\sigma(\mathcal K_T).
\]
If $c_T<1$, this spectral value is nonzero; the compact self-adjoint spectral theorem makes it an eigenvalue of finite multiplicity. A normalized eigenvector gives the stated identities. If $c_T=1$, any normalized control with a single nonzero component attains the infimum by~\eqref{eq:energy}.
\end{proof}

\subsection{Illustrations using Galerkin}\label{sec:numerics}

We approximate \eqref{eq:optimal-constant} using the nested spaces spanned by
\[
e_n(t)=\frac{\sqrt{2T}}{n\pi}\sin\frac{n\pi t}{T},
\qquad n=1,\ldots,N.
\]
These functions are orthonormal for the derivative inner product on $H_0^1(0,T)$. The resulting $MN\times MN$ energy matrix has identity diagonal blocks and off-diagonal blocks
\[
B_N(a_{j\ell}/T),\qquad
\bigl(B_N(a/T)\bigr)_{mn}=\beta_a(e_m,e_n).
\]
Let $c_{T,N}$ denote its smallest eigenvalue. The Rayleigh--Ritz principle and density imply
\begin{equation}\label{eq:galerkin-convergence}
c_T\le c_{T,N+1}\le c_{T,N},\qquad
c_{T,N}\longrightarrow c_T.
\end{equation}
To justify the last assertion, approximate any nonzero $q\in\U$ in the derivative norm by its sine partial sums. Continuity of the energy form implies convergence of their Rayleigh quotients to that of $q$. Choosing $q$ with quotient arbitrarily close to $c_T$ proves that the decreasing Galerkin limit is no larger than $c_T$; the reverse inequality follows from restriction of the infimum.

The matrix entries are evaluated by closed-form integrals. For $\alpha=a/T\in(0,1)$,
\begin{align}\label{eq:galerkin-entry}
\bigl(B_N(\alpha)\bigr)_{mn}
={}&\frac1{\alpha n\pi}\int_0^{1-\alpha}
\bigl[\cos(m\pi x)\sin(n\pi(x+\alpha))\notag\\
&\hspace{39mm}-\cos(m\pi(x+\alpha))\sin(n\pi x)\bigr]\dd x.
\end{align}
For $\alpha\ge1$, $B_N(\alpha)=0$. Product-to-sum identities reduce \eqref{eq:galerkin-entry} to the primitives
\[
\int_0^L\cos(k\pi x)\dd x=L\sinc(k\pi L),\qquad
\int_0^L\sin(k\pi x)\dd x=
\begin{cases}
\dfrac{1-\cos(k\pi L)}{k\pi},& k\ne0,\\[1mm]
0,&k=0.
\end{cases}
\]
The numerator $1-\cos(k\pi L)$ is implemented as $2\sin^2(k\pi L/2)$ to avoid subtractive cancellation. Reflection about $T/2$ gives zero entries whenever $m$ and $n$ have different parity.

For an equidistant configuration, the spatial matrix $J-I$ has eigenvalues $M-1$ and $-1$. Hence the Galerkin value can be computed from the extreme eigenvalues $\mu_{\min},\mu_{\max}$ of $B_N(a/T)$:
\[
c_{T,N}=\min\{1+(M-1)\mu_{\min},\,1-\mu_{\max}\}.
\]
For the collinear triple, reflection in the middle source reduces the matrix to an $N\times N$ block $I-B_N(2a/T)$ and a $2N\times2N$ block
\[
\begin{pmatrix}
I+B_N(2a/T)&\sqrt2B_N(a/T)\\
\sqrt2B_N(a/T)&I
\end{pmatrix}.
\]
These reductions preserve all eigenvalues relevant to the minimum.

\begin{figure}[tbp]
\centering
\includegraphics[width=\textwidth]{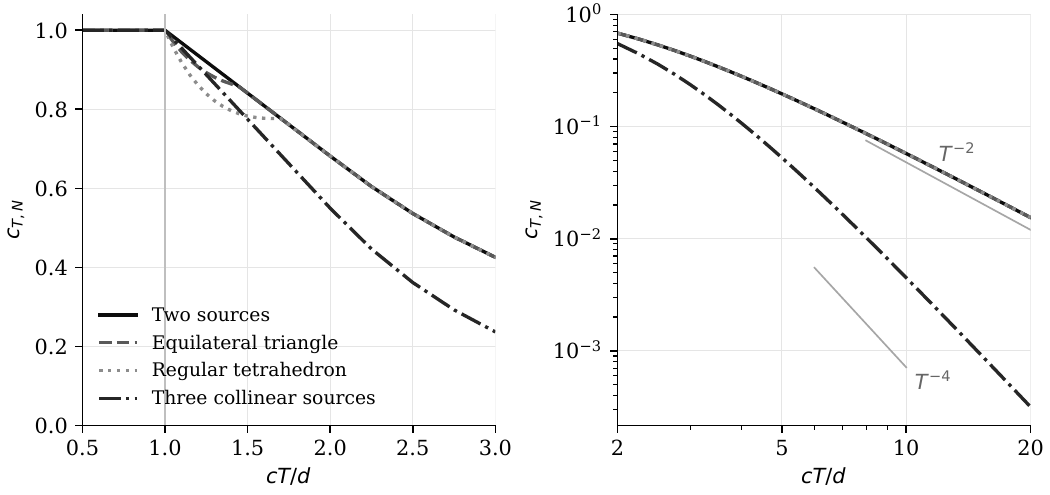}
\caption{Galerkin approximations $c_{T,N}$ with $N=160$ and $d$ the minimum source separation. Left: the short-time transition for a pair, an equilateral triangle, a regular tetrahedron, and a collinear triple with spacing $d$. Right: the same configurations on logarithmic axes. The two thin reference segments indicate slopes $T^{-2}$ and $T^{-4}$; they are not bounds. The different exponents are proved in Corollary~\ref{cor:simplex-rate} and Theorem~\ref{thm:collinear-rate}.}
\label{fig:coercivity}
\end{figure}

\begin{table}[tbp]
\centering\small
% Inlined Galerkin convergence table
\begin{tabular}{lrrrrr}
\hline
Configuration & $cT/d$ & $N=40$ & $N=80$ & $N=160$ & $N=320$ \\
\hline
Two sources & 2 & 0.6816901 & 0.6816901 & 0.6816901 & 0.6816901 \\
Two sources & 10 & 0.0575299 & 0.0575252 & 0.0575246 & 0.0575245 \\
Two sources & 20 & 0.0154432 & 0.0154379 & 0.0154372 & 0.0154371 \\
Tetrahedron & 1.5 & 0.7829119 & 0.7829052 & 0.7829044 & 0.7829043 \\
Collinear triple & 10 & 0.0044883 & 0.0044862 & 0.0044860 & 0.0044859 \\
Collinear triple & 20 & 0.0003170 & 0.0003163 & 0.0003163 & 0.0003162 \\
\hline
\end{tabular}

% End inlined Galerkin convergence table
\caption{Nested-space convergence of the computed Galerkin minima. Entries are rounded to seven decimal places.}
\label{tab:galerkin-convergence}
\end{table}

%The script \texttt{numerics/reproduce.py} generates Figure~\ref{fig:coercivity}, Table~\ref{tab:galerkin-convergence}, and the accompanying CSV data. It uses $N=160$ for the curves and $N=40,80,160,320$ for the convergence table. Selected closed-form entries are independently compared with adaptive quadrature, and the short-time identity and nested-space monotonicity are checked. The inequalities in \eqref{eq:galerkin-convergence} refer to exact Galerkin matrices. The supplied computations use ordinary floating-point arithmetic; their rounded values are not certified interval bounds. The numerical illustrations neither replace the coercivity proof nor establish a sharp long-time coefficient.
% === END QUANTITATIVE ===

% === BEGIN CONCLUSION (from supplied modular source) ===
\section{Conclusions}\label{sec:conclusion}\label{sec:open}

\noindent The finite-duration point-source propagator has a closed, explicitly characterizable range in terminal energy space. Its inverse is unique and stable once temporal support and endpoint conditions are imposed. The source geometry enters through an angular Gram matrix, while the duration of actuation determines how closely amplitude histories can concentrate near its zero-frequency degeneracy. The energy identity, radiation characterization and reconstruction formulas therefore describe complementary aspects of the same terminal-state synthesis problem.

\noindent The present analysis concerns the ideal control level, in which the amplitudes are independently prescribed. The physical control level for resonant bubbles, where these amplitudes must be generated by admissible incident waves, is an important further problem that will be addressed in forthcoming work.
\bigskip

\noindent\textbf{AI Declaration}\par

\noindent AI tools were used to assist with language editing, mathematical review, and the presentation/optimization of selected proofs. The authors are responsible for all mathematical results, references, and the final version of the manuscript.
% === END CONCLUSION ===

\raggedbottom
% === BEGIN EMBEDDED BIBLIOGRAPHY ===

% === END EMBEDDED BIBLIOGRAPHY ===
\end{document}